\RequirePackage{fix-cm}
\documentclass[smallextended]{svjour3}       
\smartqed  
\usepackage{amsmath, amssymb, amsfonts}
\usepackage{bm}
\usepackage{enumerate}
\usepackage{hyperref} 
\hypersetup{
	colorlinks=true,
	citecolor=blue,
	linkcolor=red}
\newtheorem{thm}{Theorem}
\newtheorem{prop}{Proposition}

\newtheorem{lm}{Lemma}

\newtheorem{defi}{Definition}
\numberwithin{equation}{section}

\newcommand{\hf}{\frac{1}{2}}
\newcommand{\nabh}{\nabla_{\! h}}
\newcommand{\Dh}{\Delta_{\! h}}
\newcommand{\bhu}{\hat{\mbox{\textbf{u}}}}
\newcommand{\bv}{\mathbf{v}}
\newcommand{\eipC}[2]{\left< #1 , #2 \right>_{c}}
\newcommand{\eipew}[2]{\left< #1 , #2 \right>_{ew}}
\newcommand{\eipns}[2]{\left< #1 , #2 \right>_{ns}}
\newcommand{\eipvec}[2]{\left< #1 , #2 \right>}
\newcommand{\nrminf}[1]{\left\| #1 \right\|_{\infty}}
\newcommand{\nrm}[1]{\left\| #1 \right\|}
\newcommand{\n}{\mbox{\boldmath$n$}}

\def\0{\mbox{\boldmath $0$}}
\newcommand\dt {{\Delta t}}
\newcommand{\bu}{\mathbf{u}}
\newcommand{\bn}{\mathbf{n}}
\newcommand{\bU}{\mathbf{U}}
\def\u{\mbox{\boldmath $u$}}    
\newcommand{\eu}{{\bf e}_{\textbf{u}}}
\newcommand{\ehu}{{\hat{\textbf{e}}}_{\textbf{u}}}
\newcommand{\ephi}{e_{\phi}}

\allowdisplaybreaks[2]
\begin{document}

\title{Convergence analysis of a second-order SAV-ZEC scheme for the Cahn-Hilliard-Navier-Stokes system\thanks{This work is supported in part by NSFC Youth (12101059) and the Guangdong Provincial Key Laboratory IRADS (2022B 1212010006, R 0400001-22).}
}

\titlerunning{Convergence Analysis of SAV-ZEC for CHNS}        

\author{Jingwei Sun         \and
        Zeyu Xia* \and Wei Zhang 
}

\authorrunning{J. Sun, Z. Xia, and W. Zhang} 

\institute{Jingwei Sun \at
              Department of Mathematics, National University of Defense Technology, Changsha 410073, China. \\
              \email{sjw@nudt.edu.cn}           
           \and
           Zeyu Xia \at
              School of Mathematical Sciences, University of Electronic Science and Technology of China, Chengdu 611731, China. \\
              \email{zeyuxia@uestc.edu.cn}
            \and
            Wei Zhang \at
            Guangdong Provincial Key Laboratory IRADS and Department of Mathematical Sciences, BNU-HKBU United International College, Zhuhai 519087, China.
            \email{weizhang@uic.edu.cn}
            \and 
            * represents the corresponding author.
}

\date{Received: date / Accepted: date}

\maketitle

\begin{abstract}
Incorporating the scalar auxiliary variable (SAV) method and the zero energy contribution (ZEC) technique, we analyze a linear and fully decoupled numerical scheme for the Cahn-Hilliard-Naiver-Stokes (CHNS) system. More precisely, the fully discrete scheme combines the marker-and-cell (MAC) finite difference spatial approximation and BDF2 temporal discretization, as well as the Adams-Bashforth extrapolation for the nonlinear terms, based on the SAV-ZEC reformulation. A pressure correction approach is applied to decouple the Stokes equation. Only constant-coefficient Poisson-like solvers are needed in the implementation for the resulting numerical system. The numerical scheme is unconditionally stable with respect to a rewritten total energy functional, represented in terms of one auxiliary variable in the double-well potential, another auxiliary variable to balance all the nonlinear and coupled terms, the surface energy in the original phase variable, combined with the kinematic energy part. Specifically, the error estimate for the phase variable in the $\ell^{\infty}(0,T;H_h^1)\cap\ell^2(0,T;H_h^3)$ norm, the velocity variable in the $\ell^{\infty}(0,T;\ell^2)\cap\ell^2(0,T;H_h^1)$ norm, is derived with optimal convergence rates. 
\keywords{Cahn-Hilliard-Navier-Stokes system\and marker and cell mesh\and scalar auxiliary variable\and zero energy contribution\and fully decoupled numerical scheme \and convergence and error estimate.}
\subclass{65M12 \and 65M70 \and 35K35}
\end{abstract}

\section{Introduction} 
The Cahn-Hilliard-Navier-Stokes (CHNS) system \cite{lc} is a well-known incompressible and hydrodynamically coupled model. This model describes the behavior of a fluid system undergoing phase separation, in which different components of the fluid separate into distinct regions. It is widely used in materials science for phase separation of alloys and droplet formation, in chemistry for simulating liquid mixture separation, and in biology to study cell movement and nutrient transport, see  \cite{detail1,detail2,detail3,detail4,detail5} for more details. Due to the strong coupling and complexity of this system, an efficient and accurate numerical design has always been an inevitable topic. 
An explicit form of the CHNS model in a bounded domain $\Omega\subseteq\mathbb{R}^d$ ($d=2,3$) is given by 
\begin{align} 
  & 
  \partial_t \phi + \bu\cdot \nabla\phi  = \Delta \mu 
  := \Delta ( \phi^3 - \phi - \varepsilon^2 \Delta \phi )  , \label{equation-CHNS-1} 
\\
  & 
  \partial_t \bu + \bu \cdot \nabla \bu + \nabla p = \nu \Delta \bu + \lambda \mu \nabla \phi , 
  \label{equation-CHNS-2}  
\\
  & \nabla \cdot \bu = 0,  \label{equation-CHNS-3}  
\end{align}
with no-flux boundary condition for the phase variables, and no-penetration, free-slip boundary condition for the velocity vector: 
\begin{equation}
\partial_n\phi = \partial_n\mu = 0,\quad \bu\cdot\bn=\partial_n(\bu\cdot\boldsymbol{\tau})=0,\quad\mbox{on}\,\partial_n\Omega\times(0,T].
\end{equation}   
Here, $\bu$ is the advective velocity, $p$ is the pressure variable, $\nu>0$ is the viscosity parameter, $\phi\approx\pm1$ corresponds to two different fluids, $\lambda>0$ is the mixing energy density associated with surface tension, $\mu$ stands for the chemical potential, and the parameter $\varepsilon>0$ represents the interfacial thickness. Such a hydrodynamic phase field model given by \eqref{equation-CHNS-1}-\eqref{equation-CHNS-3} is thermodynamically consistent with the second law of thermodynamics \cite {thermodynamics1,thermodynamics2,thermodynamics3} and respects an energy dissipation property:
\begin{equation} 
\begin{aligned}
 & \frac{\mathrm{d}}{\mathrm{d}t} E(\phi, \bu) = - \int_{\Omega}|\nabla \mu|^2\mathrm{d}\mathbf{x} - \frac{\nu}{\lambda} \int_{\Omega}|\nabla \bu |^2\mathrm{d}\mathbf{x} \le 0 ,\\ 
 \mbox{where}~&E(\phi, \bu)  :=  \int_\Omega \Bigl( \frac14 ( \phi^2 -1 )^2  + 1 + \frac{\varepsilon^2}{2}  | \nabla \phi |^2  \Bigr) \mathrm{d} {\bf x} + \frac{1}{2 \lambda} \int_{\Omega} |\bu |^2\mathrm{d}\mathbf{x} .  
  \label{total energy-dissipation-1} 
  \end{aligned}
\end{equation} 
See the related PDE analysis of various phase-field-fluid coupled system \cite{lc,Abels,JS}, etc.

Many successful attempts have been made to design accurate, efficient, and stable numerical algorithms for the CHNS system in the past two decades, see, e.g. \cite{feng,Shen1,Shen2,HanD,Zhaoj1,Lixiaoli,Lixiaoli2,2,yi1,yi2} and the references therein. Among them, Shen and Yang \cite{Shen1,Shen2} constructed several numerical schemes, established discrete energy laws, while the error estimate was not available. Han and Wang \cite{HanD} proposed and analyzed a second order in time method, based on the convex splitting idea and pressure projection technique. Li and Shen \cite{Lixiaoli} constructed a second order weakly-coupled, linear, energy-stable SAV-MAC scheme, the corresponding convergence analysis was established for the Cahn-Hilliard-Stokes system. Yi et al. \cite{yi1,yi2} also utilized the SAV approach to develop long-time stable numerical schemes by combining the FEM spatial discretization and the Euler semi-implicit integrator to the nonlinear coupled term. In all these works, a coupled nonlinear system has to be solved at each time step, which turns out to be a complicated process in the numerical implementation. The SAV numerical schemes have been extensively studied for many gradient flow models~\cite{YangJ2018,YangJ2019}, while the numerical design and theoretical analysis of a fully decoupled SAV approach for the CHNS system remains a challenging issue. 
On the other hand, there have been some successful attempts (\cite{yang1,yang2,yang3,yang4,yang5,yang6}) that utilized the so-called zero-energy-contribution (ZEC) technique in recent years, while none of these works have provided a theoretical proof for the error analysis. The essential difficulty arises from many auxiliary variables involved in the system, explicit treatment of nonlinear terms, as well as the extra splitting error due to the decoupling of pressure from velocity. 
In particular, we would like to highlight that Li and Shen \cite{Lixiaoli2} developed up to second order in time, fully decoupled, and unconditionally energy-stable schemes using the multiple SAV approach and established an optimal convergence analysis. It was asserted to be the first work to provide an error estimate for fully decoupled linear scheme of the CHNS system. However, the convergence analysis for fully decoupled schemes based on the ZEC technique still remains unexplored, especially for second order in time schemes. Therefore, the main purpose of this work is to derive a rigorous convergence analysis of the scheme equipped with the ZEC feature.


The key point in the theoretical analysis is to perform the error estimate for the phase variable in the $\ell^\infty (0, T; H_h^1)$ norm, and the velocity vector in the $\ell^\infty(0, T; \ell^2) \cap \ell^2 (0, T; H_h^1)$ norm, following the corresponding norms in the energy estimate, as well as the square estimate for the two auxiliary variables. Also see the related convergence analysis for various Cahn-Hilliard-fluid models~\cite{cs1,cs2,Chen2024,Chen2016,Diegel2017,GuoY2024,LiuY2017,WangC2024}, etc. In addition, an intermediate velocity vector is introduced in the consistency analysis, to facilitate the theoretical derivation associated with the pressure correction approach.  A discrete $W_h^{1, \infty}$ bound is available to the exact and constructed profiles, and an \emph{a-priori} assumption is made for the numerical error functions for the velocity and phase variables, which will be recovered in the next time step. Subsequently, a mathematical induction is employed to complete the proof. Several preliminary nonlinear error estimates are derived, and we carefully obtain the corresponding error inner product bounds for each variable. These results benefit from the fully decoupled nature induced by the ZEC technique, and a combination of these results lead to the convergence estimate of the full numerical system.

The rest of the paper is organized as follows. In Section 2, the PDE system is equivalently reformulated by the SAV and ZEC approaches. In Section 3, we propose the fully discrete finite difference scheme and state the main theoretical result. The optimal rate convergence analysis and error estimates are presented in Section 4. Finally, some concluding remarks are made in Section 5.

\section{Reformulation}
First, a nonlocal variable $q(t)$ is introduced, which satisfies a special ODE system: 
\begin{equation}\label{defi-Q-1} 
\left\{\begin{aligned}
&  q' (t) = ( \bu \cdot \nabla \phi , \mu ) - ( \mu \nabla \phi , \bu ) 
 + \lambda^{-1} ( \bu \cdot \nabla \bu , \bu ) \\
& q \mid_{t=0} \equiv 1,
\end{aligned}\right.
\end{equation}
under the condition of $\nabla\cdot \bu = 0$. Here, $(\cdot , \cdot)$ denotes the standard $L^2$ inner product. Utilizing the ZEC property satisfied by the advection and surface tension terms, it is easy to see that at the continuous level, the ODE \eqref{defi-Q-1} is equivalent to $q'(t) \equiv 0$, $q \mid_{t=0} \equiv 1$, 
so that the analytic solution to~\eqref{defi-Q-1} gives $q (t) \equiv 1$ for any $t >0$. 

Second, to derive a linear numerical scheme while preserving an alternate energy dissipation, an auxiliary variable is introduced: 
\begin{equation} 
  r (t) = \sqrt{ E_1 (\phi (t)) } , \quad E_1 (\phi) = \int_\Omega \, \Big( \frac14 \phi^4 - \frac12 \phi^2 + \frac54 \Big) \, \mathrm{d} {\bf x} ,\label{defi-R-1} 
\end{equation} 
where the constant $\frac54$ guarantees the radicand is always positive. It is obvious that $E_1 (\phi) \ge | \Omega |$ for any $\phi$. Then, by combining the two auxiliary variables $q$ and $r$, the system \eqref{equation-CHNS-1}-\eqref{equation-CHNS-3} is reformulated as 
\begin{align} 
  & 
  \partial_t \phi + q \bu \cdot \nabla \phi  = \Delta \mu , \quad 
  \mu = \frac{r}{\sqrt{E_1 (\phi) }} ( \phi^3 - \phi ) - \varepsilon^2 \Delta \phi  , \label{equation-CHNS-SAV-ZEC-1} 
\\
  & 
  r_t = \frac{1}{2 \sqrt{E_1 (\phi) }}  ( \phi^3 - \phi , \phi_t ) ,  \label{equation-CHNS-SAV-ZEC-2} 
\\
  & 
  \partial_t \bu + q \bu \cdot \nabla \bu + \nabla p = \nu \Delta \bu + \lambda q \mu \nabla \phi , 
  \label{equation-CHNS-SAV-ZEC-3}  
\\
  & 
  q_t = ( \bu \cdot \nabla \phi , \mu ) - ( \mu \nabla \phi , \bu ) 
 + \lambda^{-1} ( \bu \cdot \nabla \bu , \bu ) ,   \label{equation-CHNS-SAV-ZEC-4}  
\\
  & \nabla \cdot \bu = 0 .  \label{equation-CHNS-SAV-ZEC-5}  
\end{align}
We give some detailed descriptions of the reformulated system. 
\begin{remark}
To derive the reformulated system \eqref{equation-CHNS-SAV-ZEC-1}-\eqref{equation-CHNS-SAV-ZEC-5}, modifications have been made to the original system \eqref{equation-CHNS-1}-\eqref{equation-CHNS-3}. First, we rewrite \eqref{equation-CHNS-1} using a new auxiliary variable $r$ and take its time derivative, leading to \eqref{equation-CHNS-SAV-ZEC-2}. Second, we incorporate the ODE \eqref{defi-Q-1} into the system. Under the divergence-free and boundary conditions of the system, \eqref{defi-Q-1} is equivalent to $q \equiv 1$. To ensure the advection and surface tension terms satisfy the ZEC property, we multiply these terms by $q$. Consequently, the newly obtained PDE system, formulated in terms of the variables $(\bu, p, \phi, \mu, q, r)$, is equivalent to the original PDE system \eqref{equation-CHNS-1}-\eqref{equation-CHNS-3} formulated in $(\bu, p, \phi, \mu)$. Additionally, introducing the auxiliary variable $r$ simplifies the complex nonlinear terms in the chemical potential to a remarkably simple form, as achieved by the SAV method. This approach effectively linearizes the nonlinear terms, as demonstrated in prior studies (see, e.g., \cite{sav1,sav2,sav7}).
\end{remark}
\begin{remark}
The introduction of $q$ aims to decouple the computation of variables in the system, thereby enhancing the flexibility of algorithm design. This approach allows the nonlinear terms to be treated independently without requiring a strict cancellation condition between them. Such a decoupling is particularly beneficial to achieve stability and efficiency in the numerical design.
\end{remark}
\section{The numerical scheme and main theoretical result} 
\subsection{The finite difference spatial discretization}
For simplicity, we only consider the two dimensional domain $\Omega=(0,1)^2$. The three dimensional case could be similarly extended. In this domain, we denote the uniform spatial grid size $h=\frac{1}{N}$, with $N$ a positive integer. To facilitate the theoretical analysis, the marker and cell (MAC) grid 
is used: the phase variable $\phi$, the chemical potential $\mu$ and the pressure field $p$ are defined on the cell-centered mesh points $\left(\left(i+\hf \right)h,\,\left(j+\hf \right) h \right),\ 0\leq i,\ j \le N$; for the velocity field $\bu = (u^x,\,u^y)$, the $x$-component of the velocity is defined at the east-west cell edge points $\left(ih,\,\left(j+\hf \right) h \right),\ 0\leq i\leq N+1,\ 0\leq j \le N$, 
while the $y$-component of the velocity is located at the north-south cell edge points $\left(\left(i+\hf \right)h,\,jh\right)$. 

For a function $f(x,y)$, the notation $f_{i+\hf ,\, j+\hf }$ represents the value of $f((i+\hf )h,$ $\left(j+\hf \right)h)$. Of course, $f_{i+\hf ,\, j}$, $f_{i,\, j+\hf }$ could be similarly introduced. In turn, the following difference operators are introduced:
\begin{align}
	&(D^c_xf)_{i,\, j+\hf } = \frac{f_{i+\hf ,\, j+\hf }-f_{i-\hf ,\, j+\hf }}{h},\quad
	(D^c_yf)_{i+\hf ,\, j} = \frac{f_{i+\hf ,\, j+\hf }-f_{i+\hf ,\, j-\hf }}{h},	\label{center-diff-operator}\\
	&(D^{ew}_xf)_{i+\hf ,\, j+\hf } = \frac{f_{i+1,\, j+\hf }-f_{i,\, j+\hf }}{h},\quad
	(D^{ew}_yf)_{i,\, j} = \frac{f_{i,\, j+\hf }-f_{i,\, j-\hf }}{h},	\label{ew-diff-operator}\\
	&(D^{ns}_xf)_{i,\, j} = \frac{f_{i+\hf ,\, j}-f_{i-\hf ,\, j}}{h},\qquad\quad\quad
	(D^{ns}_yf)_{i+\hf ,\, j+\hf } = \frac{f_{i+\hf ,\, j+1}-f_{i+\hf ,\, j}}{h}.	\label{ns-diff-operator}
\end{align}
The boundary formulas may vary with different boundary conditions. 
With homogeneous Neumann boundary condition, \eqref{center-diff-operator} becomes
\begin{equation} \label{boundary-formula-1} 
	(D^c_xf)_{0,\, j+\hf } = (D^c_xf)_{N,\, j+\hf } = (D^c_yf)_{i+\hf ,\, 0} = (D^c_yf)_{i+\hf ,\, N} = 0.
\end{equation} 
The associated formulas for \eqref{ew-diff-operator}-\eqref{ns-diff-operator} could be analogously derived. 

In turn, with a careful evaluation of boundary differentiation formula~\eqref{boundary-formula-1}, the discrete boundary condition associated with cell-centered function is given by the following definition, in which the ``ghost" points are involved. The boundary formulas for the edge-centered function could be similarly derived.

	\begin{defi} \label{defi: BC} 
A cell-centered function $\phi$ is said to satisfy homogeneous Neumann boundary condition, and we write $\boldsymbol{n} \cdot \nabla_{h} \phi=0$, iff $\phi$ satisfies
	\begin{equation*} 
	\begin{aligned}
 &\phi_{-\frac12,j+\frac12}=\phi_{\frac12,j+\frac12}, \qquad\phi_{N+\frac12,j+\frac12}=\phi_{N-\frac12,j+\frac12} ,\\  
&\phi_{i+\frac12,-\frac12}=\phi_{i+\frac12,\frac12} , \qquad~\phi_{i+\frac12,N+\frac12}=\phi_{i+\frac12,N-\frac12} .  
\end{aligned}
	\end{equation*} 
A discrete function $\boldsymbol{f}=(f^{x}, f^{y})^{T}$, with two components evaluated at east-west and north-south mesh points, is said to satisfy no-penetration boundary condition, $\boldsymbol{n}\cdot\boldsymbol{f}=0$, iff we have
	\begin{equation*} 
 f_{0, j+\frac12}^x= f_{N, j+\frac12}^x=0, \qquad  
 f_{i+\frac12, 0}^{y} = f_{i+\frac12, N}^{y}=0 ,  
	\end{equation*} 
and it is said to satisfy free-slip boundary condition iff we have
	\begin{equation*} 
 f_{i, -\frac12}^{x}=f_{i, \frac12}^{x},  \qquad 
 f_{i, N+\frac12}^{x}=f_{i, N-\frac12}^{x}, \qquad 
 f_{-\frac12, j}^{y}=f_{\frac12, j}^{y} , \qquad
 f_{N+\frac12, j}^{y}=f_{N-\frac12, j}^{y}. 
	\end{equation*} 
\end{defi} 

In addition, the long stencil difference operator is also defined on the east-west cell edge points and north-south cell edge points:
\begin{align}\label{long-stencil-diff-operator}
	(\tilde{D}_xf)_{i,\, j+\hf } = \frac{f_{i+1,\, j+\hf }-f_{i-1,\, j+\hf }}{2h},\qquad
	(\tilde{D}_yf)_{i+\hf ,\, j} = \frac{f_{i+\hf ,\, j+1}-f_{i+\hf ,\, j-1}}{2h}.
\end{align}
With 
homogeneous Dirichlet boundary condition, \eqref{long-stencil-diff-operator} could be written as 
\begin{align}
	(\tilde{D}_xf)_{0,\, j+\hf } &= \frac{f_{1,\, j+\hf }-f_{-1,\, j+\hf }}{2h} = \frac{f_{1,\, j+\hf }}{h},\\
	(\tilde{D}_xf)_{N,\, j+\hf } &= \frac{f_{N+1,\, j+\hf }-f_{N-1,\, j+\hf }}{2h} = -\frac{f_{N-1,\, j+\hf }}{h},\\
	(\tilde{D}_yf)_{i+\hf ,\, 0} &= \frac{f_{i+\hf ,\, 1}-f_{i+\hf ,\, -1}}{2h} = \frac{f_{i+\hf ,\, 1}}{h},\\
	(\tilde{D}_yf)_{i+\hf ,\, N} &= \frac{f_{i+\hf ,\, N+1}-f_{i+\hf ,\, N-1}}{2h} = -\frac{f_{i+\hf ,\, N-1}}{h}.
\end{align}
For a grid function $f$, the discrete gradient operator is defined as 
\begin{align}
	&\nabh f = \left( (D^{\ell}_xf ),\  (D^{\ell}_yf )\right)^T,
\end{align}
where $\ell = c,\ ew,\ ns$ may depend on the choice of $f$. The discrete divergence operator of a vector gird function $\bu$, defined on the cell-centered points, turns out to be  
\begin{equation}
	\left(\nabh\cdot \bu\right)_{i+\hf ,\, j+\hf } = 
		(D^{ew}_xu^x)_{i+\hf ,\, j+\hf } + (D^{ns}_yu^y)_{i+\hf ,\, j+\hf }.
\end{equation}
The five point standard Laplacian operator is straightforward: 
\begin{equation}
	(\Dh f)_{r,\, s} = \frac{f_{r+1,\,s}+f_{r-1,\,s}+f_{r,\,s+1}+f_{r,\,s-1}-4f_{r,\,s}}{h^2},
\end{equation}
where $(r,\,s)$ may refer to $(i+\hf ,\,j+\hf )$, $(i+\hf ,\,j)$ and $(i,\,j+\hf )$.

For $\bu = (u^x,\,u^y)^T$, $\bv = (v^x,\,v^y)^T$, located at the staggered mesh points respectively, and the cell centered variables $\phi$, $\mu$, the nonlinear terms are evaluated as follows 
\begin{eqnarray} 
	\bu \cdot \nabla_h \bv &&= \left( \begin{array}{c} 
		u^x_{i,\,j+\hf} \tilde{D}_x v^x_{i,\,j+\hf} + {\cal A}_{xy} u^y_{i,\,j+\hf} \tilde{D}_y v^x_{i,\,j+\hf} \\ 
		{\cal A}_{xy} u^x_{i+\hf,\,j} \tilde{D}_x v^y_{i+\hf,\,j} + u^y_{i,\,j+\hf} \tilde{D}_y v^y_{i+\hf,\,j}  
	\end{array} \right)  ,  \label{FD-u-1} 
	\\
	 \mu \nabla_h \phi   &&=  \left( \begin{array}{c} 
		( D^c_x \phi \cdot {\cal A}_x \mu )_{i,\,j+\hf}   \\ 
		( D^c_y \phi \cdot {\cal A}_y \mu )_{i+\hf,\,j} 
	\end{array} \right)  ,   \label{FD-u-3} 
	\\
	\nabla_h \cdot ( \phi \bu) &&= 
	D^{ew}_x (u^x {\cal A}_x \phi)_{i+\hf,\,j+\hf} + D^{ns}_y (u^y {\cal A}_y \phi)_{i+\hf,\,j+\hf}  ,   
	\label{FD-phi-1} 
\end{eqnarray} 
where the averaging operators are given by ${\cal A}_{xy} u^x_{i+\hf,\,j} = \frac14 \bigl(  u^x_{i,\,j-\hf} + u^x_{i,\,j+\hf} 
	+ u^x_{i+1,\,j-\hf} \\ + u^x_{i+1,\,j+\hf} \bigr)$, ${\cal A}_x \phi_{i,\,j+\hf} = \frac12 \bigl(  \phi_{i-\hf,\,j+\hf} + \phi_{i+\hf,\,j+\hf} \bigr)$. A few other average terms, such as ${\cal A}_{xy} u^y_{i,\,j+\hf}$, ${\cal A}_y \phi_{i+\hf,\,j}$,  could be similarly defined.

In addition, the discrete inner product needs to be defined. Let $f$, $g$ be two grid functions evaluated on the cell-center points, the discrete $\ell^2$ inner product is given by
\begin{equation}\label{center-l2-product}
	\eipC{f}{g} = h^2\sum_{i=1}^{N}\sum_{j=1}^{N}f_{i+\hf,\,j+\hf}g_{i+\hf,j+\hf}.
\end{equation}
If $f$, $g$ are evaluated on the east-west points, \eqref{center-l2-product} becomes 
\begin{equation}\label{eastwest-l2-product}
	\eipew{f}{g} = h^2\sum_{i=1}^{N}\sum_{j=1}^{N}f_{i,\,j+\hf}g_{i,\,j+\hf}.
\end{equation}
If $f$, $g$ are evaluated on the north-south points, \eqref{center-l2-product} shifts into 
\begin{equation}\label{northsouth-l2-product}
	\eipns{f}{g} = h^2\sum_{i=1}^{N}\sum_{j=1}^{N}f_{i+\hf,\,j}g_{i+\hf,\,j}.
\end{equation}
Similarly, for two vector grid functions $\bu=\left(u^x,\,u^y\right)^T$, $\bv=\left(v^x,\,v^y\right)^T$ whose components are evaluated on east-west and north-south respectively, the vector inner product is defined as 
\begin{equation}
	\eipvec{\bu}{\bv}_1 = \eipew{u^x}{v^x} + \eipns{u^y}{v^y}.
\end{equation}
Consequently, the discrete $\ell^2$ norms, $\| \cdot \|_2$ can be naturally introduced. Furthermore, the discrete $\ell^p$, $1\leq p\leq\infty$ norms are needed in the nonlinear analysis. For $(r,\,s)=(i+\hf ,\,j+\hf )$, $(i+\hf ,\,j)$ or $(i,\,j+\hf )$, we introduce
\begin{equation}
	\nrminf{f}:=\max\limits_{r,\,s}\left|f_{r,\,s}\right|,\qquad
	\nrm{f}_p:= \Big( h^2\sum_{r=0}^{N}\sum_{s=0}^{N}\left|f_{r,\,s}\right|^p \Big)^{\frac{1}{p}},\quad1\leq p<\infty.
\end{equation} 

The discrete average is defined as $\overline{f} := \langle f , 1 \rangle_c$, for any cell centered function $f$. Moreover, an $\langle \cdot , \cdot \rangle_{-1,h}$ inner product and $\| \cdot \|_{-1, h}$ norm need to be introduced to facilitate the analysis in later sections. For any $\varphi  \in \mathring{\mathcal C}_{\Omega} := \left\{ f \middle| \ \langle f , 1 \rangle_c = 0 \right\}$, we define 
	\begin{equation} 
\langle \varphi_1, \varphi_2  \rangle_{-1,h} = \langle \varphi_1 ,  (-\Delta_h )^{-1} \varphi_2 \rangle_c , \qquad \| \varphi  \|_{-1, h } = \sqrt{ \langle \varphi ,  ( - \Delta_h )^{-1} (\varphi) \rangle_c } ,
	\end{equation} 
where the operator $\Delta_h$ is paired with discrete homogeneous Neumann boundary condition.

The following summation by parts formulas have been derived in the existing literature~\cite{Chen2024}. 

\begin{lm}  \cite{Chen2024} 
	Given two discrete grid vector functions $\bu= (u^x,\,u^y)$, $\bv=(v^x,\,v^y)$, where $u^x$, $u^y$ and $v^x$, $v^y$ are defined on east-west and north-south respectively, and two cell centered functions $f$, $g$, the following identities are valid, if $\bu$, $\bv$, $f$, $g$ are equipped with periodic boundary condition, or $\bu$, $\bv$ are implemented with homogeneous Dirichlet boundary condition and homogeneous Neumann boundary condition is imposed for $f$ and $g$:
	\begin{align} 
		\eipvec{\bu}{\nabh f} _1&= 0,\quad \mbox{if}~~\nabh\cdot \bu=0,  
		 \label{summation-2} \\
		-\eipvec{\bv}{\Dh\bv}_1 &=\nrm{\nabh\bv}_2^2, \quad 
		\eipC{f}{\Dh f} =\nrm{\nabh f}_2^2,  \label{summation-3} \\
		-\eipC{g}{\nabh\cdot\left( f\bu\right)}
		 &=\eipvec{\bu}{ f \nabh g}_1 .  \label{summation-5} 
	\end{align}
\end{lm}

The following Poincar\'{e}-type inequality and discrete Sobolev interpolation inequality will be useful in the later analysis. 

	\begin{prop} \label{prop: Poincare} 
(1) There are constants $C_0 >0$, $C_1 > 0$, independent of $h>0$, such that $\nrm{\phi}_2 \le C_0 \nrm{\nabla_h\phi}_2$, 
for all $\phi\in \mathring{\mathcal C}_{\Omega} := \left\{ f \middle| \ \langle f , 1 \rangle_c = 0 \right\}$. Moreover, we have a discrete Sobolev interpolation inequality: 
\begin{equation} 
  \| \phi \|_\infty  \le C_1 \| \phi \|_2^\frac23 \cdot \| \nabla_h \Delta_h \phi \|_2^\frac13 , \qquad 
  \| \nabla_h \phi \|_\infty \le C_1 \| \nabla_h \phi \|_2^\frac12 \cdot \| \nabla_h \Delta_h \phi \|_2^\frac12 . 
  \label{Sobolev-0} 
\end{equation} 

(2) For a velocity vector $\bv$, with a discrete no-penetration boundary condition $\bv \cdot \n =0$ on $\partial \Omega$, a similar Poincar\'e inequality is also valid: $\nrm{\bv}_2 \le C_0 \nrm{\nabla_h\bv}_2$, with $C_0$ only dependent on $\Omega$. In addition, the following discrete Sobolev interpolation inequality is available: 
\begin{equation} 
  \| \bv \|_4 \le C_2 \| \bv \|_2^\frac12 \cdot \| \nabla_h \bv \|_2^\frac12 , \quad \mbox{with $C_2$ only dependent on $\Omega$}. 
  \label{Sobolev-1} 
\end{equation}  
	\end{prop} 
	
In fact, the proof of~\eqref{Sobolev-1} has been presented in an existing work~\cite{Guan2014a}, under the periodic boundary condition. An extension to the case of no-penetration boundary condition would be straightforward, and the technical details are skipped for the sake of brevity. 

A discrete version of $E_1 (\phi)$ is defined as 
\begin{equation} 
  E_{1, h} (\phi) = \langle \frac14 \phi^4 - \frac12 \phi^2 + \frac54 , 1 \rangle_c . 
  \label{defi-E1-1} 
\end{equation} 
Again, it is clear that $E_{1, h} (\phi) \ge | \Omega| $ for any $\phi$. 
\subsection{The second order accurate numerical scheme}
A second order accurate fully discrete numerical scheme is proposed as follows. 
\begin{align} 
  & 
  \frac{\frac32 \phi^{n+1} - 2 \phi^n + \frac12 \phi^{n-1}}{\dt} + q^{n+1} \nabla_h \cdot ( \phi^* \bu^* ) 
  = \Delta_h \tilde{\mu}^{n+1} , \label{scheme-CHNS-SAV-ZEC-1} 
\\
  & 
  \tilde{\mu}^{n+1} = \frac{r^{n+1}}{\sqrt{E_{1, h} (\phi^*) }} ( (\phi^*)^3 - \phi^* ) 
   - \varepsilon^2 \Delta_h \phi^{n+1}  , \label{scheme-CHNS-SAV-ZEC-2} 
\\
  & 
  \frac{\frac32 r^{n+1} - 2 r^n + \frac12 r^{n-1}}{\dt} = \frac{1}{2 \sqrt{E_{1,h} (\phi^*) }}  
  \Big\langle ( \phi^* )^3 - \phi^* , \frac{\frac32 \phi^{n+1} - 2 \phi^n + \frac12 \phi^{n-1}}{\dt} \Big\rangle_c ,  
  \label{scheme-CHNS-SAV-ZEC-3} 
\\
  & 
  \frac{\frac32 \bhu^{n+1} - 2 \bu^n + \frac12 \bu^{n-1}}{\dt}  + q^{n+1} \bu^* \cdot \nabla_h \bu^* 
  + \nabla_h p^n \label{scheme-CHNS-SAV-ZEC-4}  \\ 
  & = \nu \Delta_h \bhu^{n+1} + \lambda q^{n+1}  \tilde{\mu}^* \nabla_h \phi^*  ,  \nonumber 
\\
  & 
  \frac{\frac32 q^{n+1} - 2 q^n + \frac12 q^{n-1}}{\dt}  \label{scheme-CHNS-SAV-ZEC-5}  \\
  &= \langle \nabla_h \cdot ( \phi^* \bu^*) , \tilde{\mu}^{n+1} \rangle_c  
  - \langle \tilde{\mu}^* \nabla_h \phi^* , \bhu^{n+1} \rangle_1 
 + \lambda^{-1} \langle \bu^* \cdot \nabla_h \bu^* , \bhu^{n+1} \rangle_1 ,\nonumber  
\\
  & 
  \frac{\bu^{n+1}-\bhu^{n+1}}{\tau}+ \frac23 \nabla_h  ( p^{n+1} -p^n )=0,	
    \label{scheme-CHNS-SAV-ZEC-6}
\\ 	
  & 
  \nabla_h \cdot \bu^{n+1}=0 ,     \label{scheme-CHNS-SAV-ZEC-7}  
\end{align}
in which 
\begin{equation}
	\begin{aligned}
		& \phi^* := 2 \phi^n- \phi^{n-1}, \quad 
	\bu^* := 2 \bu^n- \bu^{n-1}, \quad 
	\tilde{\mu}^* = ( \phi^* )^3 - \phi^* - \varepsilon^2 \Delta_h \phi^* , 
	\end{aligned}
\end{equation}
with the discrete boundary conditions:
\begin{equation}\label{physical-boundary-condition}
	( \bu^{n+1} \cdot \n ) |_{\Gamma}=0,\quad \partial_n ( \bu^{n+1} \cdot {\bf \tau} ) = 0 , \quad \partial_n\phi^{n+1}|_{\Gamma}=\partial_n\mu^{n+1 }|_{\Gamma}=0.
\end{equation}

In terms of the numerical implementation, a careful calculation reveals that \eqref{scheme-CHNS-SAV-ZEC-1}-\eqref{scheme-CHNS-SAV-ZEC-5} forms a closed numerical system for $(\phi^{n+1}, \bhu^{n+1}, r^{n+1}, q^{n+1})$. More importantly, an FFT-based fast solver could be effectively applied to this numerical system. This fact has greatly improved the numerical efficiency. In addition, the projection stage \eqref{scheme-CHNS-SAV-ZEC-6}-\eqref{scheme-CHNS-SAV-ZEC-7} corresponds to a standard Poisson equation for $p^{n+1}$, with a given $\bhu^{n+1}$, which could also be implemented by an FFT-based solver. In turn, the unique solvability analysis of the combined numerical system \eqref{scheme-CHNS-SAV-ZEC-1}-\eqref{scheme-CHNS-SAV-ZEC-7} becomes straightforward.

A modified energy stability analysis could be derived, following similar ideas as in the existing works~\cite{yang1,yang2,yang3,yang4,yang5,yang6}. 

\begin{thm} \label{thm: energy stab}
	For the proposed numerical scheme \eqref{scheme-CHNS-SAV-ZEC-1}-\eqref{scheme-CHNS-SAV-ZEC-7}, the following inequality holds for all $n>0$:
	\begin{equation}\label{discrete-energy-stab} 
	\begin{aligned} 
	    & 
	    \tilde{E}_h(\phi^{n+1}, \phi^n, r^{n+1} , r^n , q^{n+1} , q^n , \bu^{n+1}, p^{n+1}) \\ 
		\le &\tilde{E}_h(\phi^n , \phi^{n-1}, r^n , r^{n-1} , q^n , q^{n-1} , \bu^n , p^n ) ,\quad  \mbox{where} 
\\
	  & 
			\tilde{E}_h(\phi^{n+1}, \phi^n, r^{n+1} , r^n , q^{n+1} , q^n , \bu^{n+1}, p^{n+1}) 
\\
   = & 
			\frac{\varepsilon^2}{4} ( \| \nabla_h \phi^{n+1} \|_2^2 
			+ \| \nabla_h ( 2  \phi^{n+1} - \phi^n ) \|_2^2 ) 
 + \frac12 ( | r^{n+1} |^2 + | 2 r^{n+1} - r^n |^2 ) \\
  +& \frac14 ( | q^{n+1} |^2 + | 2 q^{n+1} - q^n |^2 )  
   + \frac{1}{2 \lambda} \| \bu^{n+1} \|_2^2 
			 +\frac{\dt^2}{3 \lambda} \|\nabla_h p^{n+1} \|_2^2 . 
		\end{aligned}
	\end{equation}
\end{thm}

As a direct consequence of this energy estimate, the following functional bounds for the numerical solution becomes available: 
\begin{equation} 
\begin{aligned} 
  & 
  \frac{\varepsilon^2}{4} \| \nabla_h \phi^{n+1} \|_2^2 
 + \frac12  | r^{n+1} |^2 + \frac14 | q^{n+1} |^2  
 + \frac{1}{2 \lambda} \| \bu^{n+1} \|_2^2  
\\
 \le & 
 \tilde{E}_h(\phi^{n+1}, \phi^n, r^{n+1} , r^n , q^{n+1} , q^n , \bu^{n+1}, p^{n+1})  \\
 \le &\tilde{E}_h(\phi^0, \phi^{-1}, r^0 , r^{-1} , q^0 , q^{-1} , \bu^0, p^0)  := \tilde{C}_0 , ~\mbox{so that}~
   \| \nabla_h \phi^{n+1} \|_2 \le 2 \tilde{C}_0^\frac12 \varepsilon^{-1} , 
\\
  & 
   | r^{n+1} | \le \sqrt{2} \tilde{C}_0^\frac12 , \, \, 
   | q^{n+1} | \le 2 \tilde{C}_0^\frac12 , \, \, 
   \| \bu^{n+1} \|_2 \le \sqrt{2 \lambda} \tilde{C}_0^\frac12 , 
\end{aligned} 
  \label{energy stab-1} 
\end{equation} 
for any $n \ge 0$. 
\subsection{Preliminaries and the main theorem}
Now we proceed into the convergence analysis. For the exact solution $(\phi_e, \bu_e, p_e)$ to the CHNS system~\eqref{equation-CHNS-SAV-ZEC-1}-\eqref{equation-CHNS-SAV-ZEC-5}, we could always assume that the exact solution has regularity of class $\mathcal{R}$, with sufficiently regular initial data: 
\begin{equation}
\phi_e, \, \bu_e, \, p_e \in \mathcal{R} := H^4 \left(0,T; C_{\rm per}(\Omega)\right) \cap H^3 \left(0,T; C^2_{\rm per}(\Omega)\right) \cap L^\infty \left(0,T; C^6_{\rm per}(\Omega)\right).
	\label{assumption:regularity.1}
	\end{equation}
Meanwhile, we define $\Phi_N (\, \cdot \, ,t) := {\cal P}_N \phi_e (\, \cdot \, ,t)$, the (spatial) Fourier projection of the exact solution into ${\cal B}^K$, the space of trigonometric polynomials of degree up to and including $K$ (with $N=2K +1$), only in the Cosine wave mode in both the $x$ and $y$ directions, due to the homogeneous Neumann boundary condition.  The following projection approximation is standard: if $\phi_e \in L^\infty(0,T;H^\ell_{\rm per}(\Omega))$, for some $\ell\in\mathbb{N}$,
	\begin{equation}
\nrm{\Phi_N - \phi_e}_{L^\infty(0,T;H^k)}
   \le C h^{\ell-k} \nrm{\phi_e }_{L^\infty(0,T;H^\ell)},  \quad \forall \ 0 \le k \le \ell ,
    \, \, j= 1, 2. \label{projection-est-0}
	\end{equation}
By $\Phi_N^m$, $\phi_e^m$ we denote $\Phi_N (\, \cdot \, , t_m)$ and $\phi_e (\, \cdot \, , t_m)$, respectively, with $t_m = m\cdot \dt$. Since $\Phi_N \in {\cal B}^K$, the mass conservative property is available at the discrete level:
	\begin{equation}
\overline{\Phi_N^m} = \frac{1}{|\Omega|}\int_\Omega \, \Phi_N ( \cdot, t_m) \, d {\bf x} = \frac{1}{|\Omega|}\int_\Omega \, \Phi_N ( \cdot, t_{m-1}) \, d {\bf x} = \overline{\Phi_N^{m-1}} ,  \quad \forall \ m \in\mathbb{N}.
	\label{mass conserv-1}
	\end{equation}
On the other hand, the numerical solution of the phase variable is also mass conservative at the discrete level. 
Meanwhile, we use the mass conservative projection for the initial data:  $\phi^0 = {\mathcal P}_h \Phi_N (\, \cdot \, , t=0)$, that is
	\begin{equation}
\phi^0_{i,j} := \Phi_N (p_i, p_j, t=0)  .
	\label{initial data-0}
	\end{equation}	
In turn, the error grid function for the phase variable is defined as
	\begin{equation}
e_\phi^m := \mathcal{P}_h \Phi_N^m - \phi^m ,  \quad m=0,1,2,\ldots.
	\label{error function-1}
	\end{equation}
Therefore, it follows that  $\overline{e_\phi^m} =0$. In turn, a discrete Ponicar\'e inequality becomes available for $e_\phi^m$. 

In terms of the velocity vector, 
it is observed that the exact velocity profile $\bu_e$ is not divergence-free at a discrete level, so that its discrete inner product with the pressure gradient may not vanish. To overcome this subtle difficulty, a spatial interpolation operator is needed to ensure the exact divergence-free property of the constructed velocity vector at a discrete level. Such an operator in the finite difference discretization is highly non-standard, due to the collocation point structure, and this effort has not been reported in the existing textbook literature. A pioneering idea of this approach was proposed in an existing work~\cite{GuoY2024}, 
and other related analysis works have been reported. 
In more details, the spatial interpolation operator $\bm{\mathcal{P}}_H$ is defined as follows, for any $\bu \in H^1(\Omega)$, $\nabla\cdot \bu=0$: There is an exact stream function $\psi$ so that $\bu=\nabla^{\perp} \psi$, and we define 
\begin{equation}
\bm{\mathcal{P}}_H (\bu)=\nabla_h^{\perp} \psi =( - D_y \psi , D_x \psi )^T . 
\end{equation}
Of course, this definition ensures $\nabla_h \cdot \bm{\mathcal{P}}_H (\bu)=0$ at a point-wise level, and an $O (h^2)$ truncation error is available between the continuous velocity $\bu$ and its Helmholtz interpolation, $\bm{\mathcal{P}}_H (\bu)$. 

In turn, we take $\bU = \bm{\mathcal{P}}_H (\bu_e)$, so that this constructed vector is divergence-free at a discrete level, and it is within an $O (h^2)$ approximation to the exact profile. Meanwhile, we just take $P = p_e$ for the pressure variable. Subsequently, the associated error grid functions are defined as
	\begin{equation}
{\bf e}_{\bu}^m := \mathcal{P}_h \bU^m - \bu^m = ( e_u^m , e_v^m)^T,  \, \, \, 
e_p^m := \mathcal{P}_h P^m - p^m, \quad  m = 0 ,1 , 2,\ldots .
	\label{error function-2}
	\end{equation}
  
  The following theorem is the main result of this article.

\begin{thm}
	\label{thm: convergence}
Given initial data $\phi_e (\, \cdot \, ,t=0), \, \u_e (\, \cdot \, ,t=0) \in C^6_{\rm per}(\Omega)$, suppose the exact solution for CHNS system~\eqref{equation-CHNS-SAV-ZEC-1}-\eqref{equation-CHNS-SAV-ZEC-5} is of regularity class $\mathcal{R}$. Then, provided $\dt$ and $h$ are sufficiently small, 
we have
	\begin{equation}
	\begin{aligned}
  &\varepsilon \| \nabla_h e_\phi^n \|_2 +  \| {\bf e}_{\bu}^n \|_2 
  + \varepsilon^2 \Big(  \dt   \sum_{m=1}^{n} \| \nabla_h \Delta_h e_\phi^m \|_2^2 \Big)^{\hf} 
  + \Big( \nu \dt \sum_{m=1}^{n}  \| \nabla_h \eu^m \|_2^2 \Big)^\frac12\\
   \le& C ( \dt^2+ h^2 ) ,
   \end{aligned}
	\label{convergence-0}
	\end{equation}
for all positive integers $n$, such that $t_n=n \dt \le T$, where $C>0$ is independent of $\dt$ and $h$.
	\end{thm}

\section{The convergence analysis}

\subsection{Consistency analysis and error evolutionary system} 

The following intermediate velocity vector is defined, which is needed in the later analysis: 
\begin{equation} 
   \hat{\bU}^{n+1} = \bU^{n+1} + \frac23 \dt \nabla_h ( P^{n+1} - P^n ) .  \label{consistency-hatu-1} 
\end{equation} 
In addition, we denote $R:= r_e$ and $Q:= q_e$ as the exact scalar profiles for $r$ and $q$, respectively. For the projection solution $\Phi_N$, the constructed velocity profile $\bU$ and the exact pressure variable $P$, as well as the exact scalar variables $R$ and $Q$, a careful Taylor expansion (in both time and space) gives  
\begin{align} 
  & 
  \frac{\frac32 \Phi_N^{n+1} - 2 \Phi_N^n + \frac12 \Phi_N^{n-1}}{\dt} 
  + Q^{n+1} \nabla_h \cdot ( \Phi_N^* \bU^* ) 
  = \Delta_h M^{n+1} + \mathbf{G}_\phi^{n+1}  , \label{consistency-CHNS-SAV-ZEC-1} 
\\
  & 
  M^{n+1} = \frac{R^{n+1}}{\sqrt{E_{1, h} (\Phi_N^*) }} ( (\Phi_N^*)^3 - \Phi_N^* ) 
   - \varepsilon^2 \Delta_h \Phi_N^{n+1}  , \label{consistency-CHNS-SAV-ZEC-2} 
\\
  & 
  \frac{\frac32 R^{n+1} - 2 R^n + \frac12 R^{n-1}}{\dt}   \label{consistency-CHNS-SAV-ZEC-3} \\
  =& \frac{1}{2 \sqrt{E_{1,h} (\Phi_N^*) }}  
  \Big\langle ( \Phi_N^* )^3 - \Phi_N^* , 
   \frac{\frac32 \Phi_N^{n+1} - 2 \Phi_N^n + \frac12 \Phi_N^{n-1}}{\dt} \Big\rangle_c 
   + \mathbf{G}_r^{n+1} ,  
\nonumber
\\
  & 
  \frac{\frac32 \hat{\bU}^{n+1} - 2 \bU^n + \frac12 \bU^{n-1}}{\dt}  + Q^{n+1} \bU^* \cdot \nabla_h \bU^* 
  + \nabla_h P^n \label{consistency-CHNS-SAV-ZEC-4}  \\
  =& \nu \Delta_h \hat{\bU}^{n+1} + \lambda Q^{n+1} M^* \nabla_h \Phi_N^*  
  + \mathbf{G}_{\bu}^{n+1}  , 
  \nonumber
\\
  & 
  \frac{\frac32 Q^{n+1} - 2 Q^n + \frac12 Q^{n-1}}{\dt}
  = \langle \nabla_h \cdot ( \Phi_N^* \bU^* ) , M^{n+1} \rangle_c   \label{consistency-CHNS-SAV-ZEC-5}  
\\
  -& \langle M^* \nabla_h \Phi_N^* , \hat{\bU}^{n+1} \rangle_1   
  + \lambda^{-1} \langle \bU^* \cdot \nabla_h \bU^* , \hat{\bU}^{n+1} \rangle_1 
 + \mathbf{G}_q^{n+1} ,  \nonumber
\\
  & 
  \frac{\bU^{n+1}- \hat{\bU}^{n+1}}{\dt}+ \frac23 \nabla_h  ( P^{n+1} -P^n )=0,	
    \label{consistency-CHNS-SAV-ZEC-6}
\\ 	
  & 
  \nabla_h \cdot \bU^{n+1}=0 ,     \label{consistency-CHNS-SAV-ZEC-7}  
\end{align}
in which $\|\mathbf{G}_\phi^{n+1} \|_{-1,h}$, $| \mathbf{G}_r^{n+1} |$, $\|\mathbf{G}_{\bu}^{n+1} \|_2$, $|\mathbf{G}_q^{n+1} | \le C (\dt^2 + h^2)$, and $C$ depends on the regularity of the exact solution. The star profiles are given by 
\begin{equation}
	\begin{aligned}
		& \Phi_N^* := 2 \Phi_N^n- \Phi_N^{n-1}, \, \,  
	\bU^* := 2 \bU^n- \bU^{n-1}, \, \, 
	M^* = ( \Phi_N^* )^3 - \Phi_N^* - \varepsilon^2 \Delta_h \Phi_N^* , 
	\end{aligned}
\end{equation} 

Due to the regularity of exact solution $(\Phi, \bU, P)$, its discrete $W_h^{1,\infty}$ norm will stay bounded: 
\begin{equation} 
\begin{aligned} 
  & 
  \| \Phi_N^k \|_\infty + \| \nabla_h \Phi_N^k \|_\infty \le C^\star ,~
  \| \bU^k \|_\infty + \| \nabla_h \bU^k \|_\infty \le C^\star,~
  Q^k \equiv 1 , ~ R^k \le C^\star  ,
\\
  & 
   \| \hat{\bU}^{n+1} \|_\infty \le \| \bU^{n+1} \|_\infty + \frac23 \dt \| \nabla_h ( P^{n+1} - P^n ) \|_\infty 
   \le \tilde{C}_2 := C^* + \frac12 , 
\\
  & 
  \| \hat{\bU}^{n+1} \|_2 \le \| {\bf 1 } \|_2 \cdot \| \hat{\bU}^{n+1} \|_\infty 
  \le | \Omega |^\frac12 \tilde{C}_2 ,\,
  \| \hat{\bU}^{n+1} \|_4 \le \| {\bf 1 } \|_4 \cdot \| \hat{\bU}^{n+1} \|_\infty 
  \le | \Omega |^\frac14 \tilde{C}_2 ,
\end{aligned} 
    \label{assumption:W1-infty-1}
\end{equation} 
provided that $\dt$ is sufficiently small and for all $k \ge 0 $. In particular, the following discrete $W_h^{1,\infty}$ bounds for the star profiles are also available: 
\begin{equation} 
  \| \Phi_N^* \|_\infty + \| \nabla_h \Phi_N^* \|_\infty \le 3 C^\star ,\qquad
  \| \bU^* \|_\infty + \| \nabla_h \bU^* \|_\infty \le 3 C^\star   . 
    \label{assumption:W1-infty-2}
\end{equation} 
Similarly, since $M^{n+1}$ and $M^*$ only depends on the exact solution $\Phi$, we assume a discrete $H_h^1$ and $\| \cdot \|_\infty$ bound
\begin{equation}
  \| \nabla_h M^{n+1}\|_2 , \quad 
 \| M^* \|_\infty 
 \le C^{**} ,  \label{assumption:W1-infty-mu}
\end{equation}
with $C^{**}$ a constant only dependent on the regularity of the exact solution. Moreover, due to the regularity of $\Phi$ in time, its discrete temporal derivative turns out to be bounded in the $\| \cdot \|_\infty$ norm: 
\begin{equation} 
  \Big\| \frac{\frac32 \Phi_N^{n+1} - 2 \Phi_N^n + \frac12 \Phi_N^{n-1}}{\dt} \Big\|_\infty \le C^* . 
   \label{assumption:W1-infty-3}
\end{equation}

In addition to the error functions defined in~\eqref{error function-1}, \eqref{error function-2}, the following auxiliary error functions are introduced: 
\begin{equation} 
	\begin{aligned}
		&
		\eu^*= \bU^*-\bu^* = 2 e_{\bu}^n - e_{\bu}^{n-1} ,\, \, \,  e_\phi^*= \Phi_N^*- \phi^* = 2 e_\phi^n - e_\phi^{n-1} , \, \, \,  
\\  
		& \ehu^{n+1} = \hat{\bU}^{n+1} - \hat{\bu}^{n+1} , \,\,\, e_r^k = R^k - r^k , \, \, \,  e_q^k = Q^k - q^k, \, \, \, \\
		& e_{\mu}^{n+1} = M^{n+1} -\tilde{\mu}^{n+1} , \, \, \, 
		  e_{\mu}^* = M^* - \tilde{\mu}^* .
	\end{aligned} 
	\label{error function-3} 
\end{equation} 
In turn, subtracting the numerical system \eqref{scheme-CHNS-SAV-ZEC-1}-\eqref{scheme-CHNS-SAV-ZEC-7} from the consistency estimate \eqref{consistency-CHNS-SAV-ZEC-1}-\eqref{consistency-CHNS-SAV-ZEC-7} leads to the following error evolutionary system: 
\begin{align} 
  & 
  \frac{\frac32 e_\phi^{n+1} - 2 e_\phi^n + \frac12 e_\phi^{n-1}}{\dt} 
   +  \nabla_h \cdot ( e_\phi^* \bU^* + \phi^* \eu^* ) 
   + e_q^{n+1} \nabla_h \cdot ( \phi^* \bu^* )   \label{error-CHNS-SAV-ZEC-1}  \\ 
   = & 
    \Delta_h e_\mu^{n+1} + \mathbf{G}_\phi^{n+1}  , \nonumber 
\\
  & 
  e_\mu^{n+1} = e_\mu^{n+1, (1)} + e_\mu^{n+1, (2)} + e_\mu^{n+1, (3)}
     - \varepsilon^2 \Delta_h e_\phi^{n+1}  ,~~\mbox{where} \label{error-CHNS-SAV-ZEC-2}  
\\
  & 
  e_\mu^{n+1, (1)} = \frac{e_r^{n+1}}{\sqrt{E_{1, h} (\phi^*) }} ( (\Phi_N^*)^3 - \Phi_N^* )  , \, \,  
  e_\mu^{n+1, (2)}  = \frac{r^{n+1}}{\sqrt{E_{1, h} (\phi^*) }} ( {\cal NLP}^* -1 ) e_\phi^* ,  \nonumber 
\\
  & 
  e_\mu^{n+1, (3)} = \frac{R^{n+1} \bigl( \sqrt{E_{1, h} (\phi^*)} - \sqrt{E_{1,h} (\Phi_N^*)} \bigr)}{\sqrt{E_{1, h} (\phi^*) } 
  \cdot \sqrt{E_{1, h} (\Phi_N^*) }} ( (\Phi_N^*)^3 - \Phi_N^* )  ,\nonumber \\
 & {\cal NLP}^* = (\Phi_N^* )^2 + \Phi_N^* \cdot \phi^* + (\phi^* )^2 , \nonumber 
\\
  & 
  \frac{\frac32 e_r^{n+1} - 2 e_r^n + \frac12 e_r^{n-1}}{\dt} 
  = \frac{1}{2 \sqrt{E_{1,h} (\phi^*) }}  
  \Big\langle  ( \Phi_N^* )^3 - \Phi_N^* , 
   \frac{\frac32 e_\phi^{n+1} - 2 e_\phi^n + \frac12 e_\phi^{n-1}}{\dt} \Big\rangle_c    \label{error-CHNS-SAV-ZEC-3} 
\\
  +& \frac{1}{2 \sqrt{E_{1,h} (\phi^*) }}  
  \Big\langle ( {\cal NLP}^* -1 ) e_\phi^* , 
   \frac{\frac32 \Phi_N^{n+1} - 2 \Phi_N^n + \frac12 \Phi_N^{n-1}}{\dt} \Big\rangle_c  
   + \mathbf{G}_r^{n+1} \nonumber 
\\
  +&  \frac{\sqrt{E_{1, h} (\phi^*)} - \sqrt{E_{1,h} (\Phi_N^*)}}{2 \sqrt{E_{1, h} (\phi^*) } 
  \cdot \sqrt{E_{1, h} (\Phi_N^*) }} 
  \Big\langle ( \Phi_N^* )^3 - \Phi_N^* , 
   \frac{\frac32 \Phi_N^{n+1} - 2 \Phi_N^n + \frac12 \Phi_N^{n-1}}{\dt} \Big\rangle_c  ,  
\nonumber
\\
  & 
  \frac{\frac32 \ehu^{n+1} - 2 \eu^n + \frac12 \eu^{n-1}}{\dt}  
  + ( \bu^* \cdot \nabla_h \eu^* + \eu^* \cdot \nabla_h \bU^* ) 
  + e_q^{n+1} ( \bu^* \cdot \nabla_h \bu^* )   \label{error-CHNS-SAV-ZEC-4}  
\\
   + & \nabla_h e_p^n - \nu \Delta_h \ehu^{n+1} 
  = \lambda ( M^* \nabla_h e_\phi^* + e_{\mu}^* \nabla_h \phi^* 
  + e_q^{n+1} \tilde{\mu}^* \nabla_h \phi^*  ) 
  + \mathbf{G}_{\bu}^{n+1}  , 
  \nonumber
\\
  & 
  \frac{\frac32 e_q^{n+1} - 2 e_q^n + \frac12 e_q^{n-1}}{\dt} 
  =  \langle \nabla_h \cdot ( \bU^* e_\phi^* + \eu^* \phi^* ) , M^{n+1} \rangle_c \label{error-CHNS-SAV-ZEC-5}
\\
  + & \langle \nabla_h \cdot ( \phi^* \bu^* ) , e_\mu^{n+1} \rangle_c  
  -  \bigl(  \langle \tilde{\mu}^* \nabla_h \phi^* , \ehu^{n+1} \rangle_1 
   + \langle M^* \nabla_h e_\phi^* + e_{\mu}^* \nabla_h \phi^* , \hat{\bU}^{n+1} \rangle_1 \bigr)  
    \nonumber 
\\
 +& 
  \lambda^{-1} ( \langle \bu^* \cdot \nabla_h \bu^* , \ehu^{n+1} \rangle_1 
 + \langle \bu^* \cdot \nabla_h \eu^* + \eu^* \cdot \nabla_h \bU^* , 
 \hat{\bU}^{n+1} \rangle_1 ) 
 + \mathbf{G}_q^{n+1} ,   \nonumber  
\\
  & 
  \frac{\eu^{n+1}- \ehu^{n+1}}{\dt}+ \frac23 \nabla_h  ( e_p^{n+1} - e_p^n )=0,	
    \label{error-CHNS-SAV-ZEC-6}
\\ 	
  & 
  \nabla_h \cdot \eu^{n+1}=0 ,     \label{error-CHNS-SAV-ZEC-7}  
\end{align} 
in which the fact that $Q^{n+1} \equiv 1$ has been repeatedly applied. 

\subsection{The \emph{a-priori} assumption and preliminary estimates} 

To proceed with the convergence analysis, the following \emph{a-priori} assumption is made for the numerical error functions at the previous time steps: 
\begin{equation} \label{a priori-1}
	\| \eu^k \|_4 , \, \, \, \| \ephi^k \|_\infty + \| \nabla_h \ephi^k \|_\infty , \, \, \, \| \Delta_h \ephi^k \|_2 \leq \dt^\frac14 + h^\frac14 ,  \, \, \, k= n, n-1 . 
\end{equation} 
Such an \emph{a-priori} assumption will be recovered by the convergence analysis in the next time step,  which will be demonstrated later. 
In turn, the \emph{a-priori} assumption~\eqref{a priori-1} leads to a $W_h^{1, \infty}$ bound for the numerical solution for the phase variable, as well as a $\| \cdot \|_4$ bound for that of the velocity vector, for $k=n, n-1$:  
\begin{equation} 
\begin{aligned}  
&\| \phi^k \|_\infty + \| \nabla_h \phi^k \|_\infty \le \| \Phi_N^k \|_\infty + \| \nabla_h \Phi_N^k \|_\infty 
 + \| \ephi^k \|_\infty + \| \nabla_h \ephi^k \|_\infty \\
\le& C^* + \frac12 = \tilde{C}_2 ,  
\\
  &\| \bu^k \|_4 \le \| \bU^k \|_4 + \| \eu^k \|_4 \le \breve{C}_0 C^* + \frac12 \le \breve{C}_0 \tilde{C}_2 , 
\end{aligned}   
 \label{a priori-2} 
\end{equation} 
provided that $\dt$ and $h$ are sufficiently small, in which the $W_h^{1, \infty}$ assumption~\eqref{assumption:W1-infty-1} has been recalled. Notice that $\breve{C}_0$ is the H\"older inequality constant, $\| f \|_4 \le \breve{C}_0 \| f \|_\infty$, and such a constant only depends on $\Omega$. Moreover, the established bounds in~\eqref{a priori-2} implies that $\| M^* - \tilde{\mu}^* \|_2 \le C ( \dt^{\frac14} + h^\frac14 ) \frac12$. As a consequence, the corresponding bounds for the star numerical profiles are also valid:  
\begin{equation} 
\begin{aligned}  
  & 
\| \phi^* \|_\infty + \| \nabla_h \phi^* \|_\infty \le 3 \tilde{C}_2 ,  \, \, \, 
\| \tilde{\mu}^* \|_2 \le | \Omega |^\frac12 C^{**} + \frac12 := \tilde{C}_{2,2} , 
\\
  & 
  \| \bu^* \|_4 \le \breve{C}_1 \tilde{C}_2 ,   \, \, \, \breve{C}_1 = 3 \breve{C}_0 ,  
\end{aligned}   
 \label{a priori-3} 
\end{equation} 
due to the fact that $\phi^* = 2 \phi^n - \phi^{n-1}$, $\bu^* = 2 \bu^n - \bu^{n-1}$. 

A few preliminary nonlinear error estimates are stated below. The corresponding proofs are placed in Appendices~\ref{appA} and \ref{appB}.

\begin{prop}  \label{prop: nonlinear est} 
Assume the functional bounds \eqref{assumption:W1-infty-1}-\eqref{assumption:W1-infty-3} for the exact and constructed solutions, as well as the \emph{a-priori} assumption~\eqref{a priori-1}, the following estimates are valid: 
\begin{align} 
  & 
  \| {\cal NLP}^* \|_\infty + \| \nabla_h {\cal NLP}^* \|_\infty \le \tilde{C}_3 ,   \label{nonlinear est-0-1} 
\\
  & 
  \| ( \Phi_N^* )^3 - \Phi_N^* \|_\infty + \| \nabla_h ( ( \Phi_N^* )^3 - \Phi_N^* ) \|_\infty \le \tilde{C}_4;  
  \nonumber 
\\
  &  
  \| e_\mu^{n+1, (2)} \|_2 + \| \nabla_h e_\mu^{n+1, (2)} \|_2 
  \le \tilde{C}_5 ( \| e_\phi^* \|_2  + \| \nabla_h e_\phi^* \|_2) ;  \label{nonlinear est-0-2} 
\\
  & 
  | E_{1, h} (\phi^*) - E_{1,h} (\Phi_N^*) | , \quad 
  \Big| \sqrt{E_{1, h} (\phi^*)} - \sqrt{E_{1,h} (\Phi_N^*)} \Big| \le \tilde{C}_6 \| e_\phi^* \|_2 ;  
  \label{nonlinear est-0-3} 
\\
  & 
   \| e_\mu^{n+1, (1)} \|_2 + \| \nabla_h e_\mu^{n+1, (1)} \|_2 \le \tilde{C}_7 | e_r^{n+1} | ,   
   \label{nonlinear est-0-4}  
\\
  &   
  \| e_\mu^{n+1, (3)} \|_2 + \| \nabla_h e_\mu^{n+1, (3)} \|_2 \le \tilde{C}_8 \| e_\phi^* \|_2  ; 
  \nonumber 
\\
  & 
   \frac{1}{\sqrt{E_{1,h} (\phi^*) }}  
  \Big\langle ( {\cal NLP}^* -1 ) e_\phi^* , 
   \frac{\frac32 \Phi_N^{n+1} - 2 \Phi_N^n + \frac12 \Phi_N^{n-1}}{\dt} \Big\rangle_c
   \le \tilde{C}_9 \| e_\phi^* \|_2 ; \label{nonlinear est-0-5} 
\\
  & 
    \frac{\sqrt{E_{1, h} (\phi^*)} - \sqrt{E_{1,h} (\Phi_N^*)}}{\sqrt{E_{1, h} (\phi^*) } 
  \cdot \sqrt{E_{1, h} (\Phi_N^*) }} 
  \Big\langle ( \Phi_N^* )^3 - \Phi_N^* , 
   \frac{\frac32 \Phi_N^{n+1} - 2 \Phi_N^n + \frac12 \Phi_N^{n-1}}{\dt} \Big\rangle_c  
   \label{nonlinear est-0-6}  
 \\
   & 
   \le \tilde{C}_{10} \| e_\phi^* \|_2 ;  \nonumber  
\\
  & 
  \| e_\phi^* \bU^* +  \phi^* \eu^* \|_2 
  \le \tilde{C}_{11} (  \| e_\phi^* \|_2 +  \| \eu^* \|_2 ) ,  \quad 
  \| \bu^* \phi^* \|_2 \le \tilde{C}_{12} ;  \label{nonlinear est-0-7} 
\\
  & 
  \| M^* \nabla_h e_\phi^* + e_{\mu}^* \nabla_h \phi^* \|_2 
  \le \tilde{C}_{13} ( \| \nabla_h e_\phi^* \|_2 + \| e_{\mu}^* \|_2 ) ,  \label{nonlinear est-0-8}   
\\
  & 
   \| e_{\mu}^* \|_2^2 \le \tilde{C}_{14} \| e_\phi^* \|_2^2 
   + 2 \varepsilon^4 \| \Delta_h e_\phi^* \|_2^2 , \nonumber 
\end{align} 
in which $\tilde{C}_j$ ($3 \le j \le 14$) only depends on the regularity of the exact solution, the domain $\Omega$ and the initial data. 
\end{prop} 

In terms of the pressure correction stage~\eqref{error-CHNS-SAV-ZEC-6}-\eqref{error-CHNS-SAV-ZEC-7}, the following estimates will be helpful in the later analysis. 

\begin{prop}  \label{prop: projection est} 
For $\ehu^{n+1}$ and $\eu^{n+1}$ satisfying \eqref{error-CHNS-SAV-ZEC-6}-\eqref{error-CHNS-SAV-ZEC-7}, we have 
\begin{equation}  
\begin{aligned}
&  \| \ehu^{n+1} \|_2^2 = \| \eu^{n+1} \|_2^2 + \| \ehu^{n+1} - \eu^{n+1} \|_2^2 ,\\
 & \| \nabla_h \ehu^{n+1} \|_2^2 = \| \nabla_h \eu^{n+1} \|_2^2 + \| \nabla_h ( \ehu^{n+1} - \eu^{n+1} ) \|_2^2 . 
 \end{aligned}
  \label{projection est-0-1} 
\end{equation} 
\end{prop} 

\subsection{Error estimates} 

By using the \emph{a-priori} assumption~\eqref{a priori-1}, the resulting \emph{a-priori} bounds \eqref{a priori-2}-\eqref{a priori-3}, the regularity assumptions \eqref{assumption:W1-infty-1}-\eqref{assumption:W1-infty-3}, as well as the preliminary nonlinear error estimates stated in Proposition~\ref{prop: nonlinear est}, we are able to derive the convergence analysis of the SAV-ZEC numerical scheme. 

Taking a discrete inner product with~\eqref{error-CHNS-SAV-ZEC-1} by $\hat{e}_{\mu}^{n+1} := e_\mu^{n+1, (1)} - \varepsilon^2 \Delta_h e_\phi^{n+1}$ yields 
\begin{equation} 
\begin{aligned} 
  & 
  \frac{1}{\dt} \langle \frac32 e_\phi^{n+1} - 2 e_\phi^n + \frac12 e_\phi^{n-1} ,  - \varepsilon^2 \Delta_h e_\phi^{n+1} + e_\mu^{n+1, (1)} \rangle_c  
\\
  - &  \langle e_\phi^* \bU^* +  \phi^* \eu^* ,  \nabla_h \hat{e}_{\mu}^{n+1} \rangle_1  
   - e_q^{n+1} \langle \phi^* \bu^* ,  \nabla_h ( e_{\mu}^{n+1} - e_\mu^{n+1, (2)} 
   - e_\mu^{n+1, (3)} ) \rangle_1   
 \\
  + & \langle \nabla_h e_\mu^{n+1} , \nabla_h \hat{e}_{\mu}^{n+1} \rangle_1   
  = \langle  \hat{e}_{\mu}^{n+1} , \mathbf{G}_\phi^{n+1}  \rangle_c , 
\end{aligned} 
  \label{convergence-CHNS-SAV-ZEC-1} 
\end{equation} 
in which the summation by parts formulas, as well as the identity, $\hat{e}_{\mu}^{n+1} = e_{\mu}^{n+1} - e_\mu^{n+1, (2)} - e_\mu^{n+1, (3)}$, have been used. Meanwhile, the following equalities and the associated bounds could be carefully derived: 
\begin{align} 
  & 
  \langle \frac32 e_\phi^{n+1} - 2 e_\phi^n + \frac12 e_\phi^{n-1} ,  - \Delta_h e_\phi^{n+1} \rangle_c 
  = \langle \nabla_h ( \frac32 e_\phi^{n+1} - 2 e_\phi^n + \frac12 e_\phi^{n-1} ) ,  
  \nabla_h e_\phi^{n+1} \rangle_c      \label{convergence-CHNS-SAV-ZEC-2-1} 
\\
  \ge & 
  \frac14 ( \| \nabla_h e_\phi^{n+1} \|_2^2 - \| \nabla_h e_\phi^n \|_2^2 
  + \| \nabla_h ( 2 e_\phi^{n+1} - e_\phi^n )\|_2^2 
   - \| \nabla_h ( 2 e_\phi^n - e_\phi^{n-1} ) \|_2^2  ) ; 
\nonumber
\\
  & 
  \langle e_\phi^* \bU^* +  \phi^* \eu^* ,  \nabla_h \hat{e}_{\mu}^{n+1} \rangle_1
  \le \| e_\phi^* \bU^* +  \phi^* \eu^* \|_2 
  \cdot \| \nabla_h \hat{e}_{\mu}^{n+1} \|_2   \label{convergence-CHNS-SAV-ZEC-2-2} 
\\
  \le &  
    \tilde{C}_{11} ( e_\phi^* +  \eu^* ) \cdot \| \nabla_h \hat{e}_{\mu}^{n+1} \|_2 
  \le 8 \tilde{C}_{11}^2 ( \| e_\phi^* \|_2^2 +  \| \eu^* \|_2^2 ) 
  + \frac{1}{16} \| \nabla_h \hat{e}_{\mu}^{n+1} \|_2^2 ;  
\nonumber
\\
  & 
  \langle  \phi^* \bu^* ,  \nabla_h ( e_\mu^{n+1, (2)} + e_\mu^{n+1, (3)}  ) \rangle_1 
  \label{convergence-CHNS-SAV-ZEC-2-3} 
\\
  \le & \| \phi^* \bu^* \|_2 \cdot ( \|  \nabla_h e_\mu^{n+1, (2)} \|_2 
  + \| \nabla_h e_\mu^{n+1, (3)}  \|_2 )     \nonumber 
\\
  \le & 
    \tilde{C}_{12} \cdot ( \tilde{C}_5 + \tilde{C}_8 ) ( \| e_\phi^* \|_2 + \| \nabla_h e_\phi^* \|_2 )  \quad 
  \mbox{(by~\eqref{nonlinear est-0-2}, \eqref{nonlinear est-0-4}, \eqref{nonlinear est-0-7})};
\nonumber 
\\
  & 
  - e_q^{n+1} \langle  \phi^* \bu^* ,  \nabla_h ( e_\mu^{n+1, (2)} + e_\mu^{n+1, (3)} ) \rangle_1  
  \label{convergence-CHNS-SAV-ZEC-2-4} 
\\ 
  \le & \tilde{C}_{12} ( \tilde{C}_5 + \tilde{C}_8 ) | e_q^{n+1} | \cdot ( \| e_\phi^* \|_2 + \| \nabla_h e_\phi^* \|_2 )  \nonumber 
\\
  \le &  
   \tilde{C}_{12} ( \tilde{C}_5 + \tilde{C}_8 ) ( \frac12 | e_q^{n+1} |^2 
  + \| e_\phi^* \|_2 ^2 + \| \nabla_h e_\phi^* \|_2^2 );
 \nonumber 
\\
  & 
   \langle \nabla_h e_\mu^{n+1} , \nabla_h \hat{e}_{\mu}^{n+1} \rangle_c  
   =  \| \nabla_h \hat{e}_{\mu}^{n+1} \|_2^2   
    + \big\langle \nabla_h ( e_\mu^{n+1,(2)} + e_\mu^{n+1, (3)} ) , \nabla_h \hat{e}_{\mu}^{n+1} \big\rangle_c  
  \label{convergence-CHNS-SAV-ZEC-2-5} 
\\
  \ge &  
    \frac78 \| \nabla_h \hat{e}_{\mu}^{n+1} \|_2^2   
    - 2 \big\| \nabla_h ( e_\mu^{n+1,(2)} + e_\mu^{n+1, (3)} ) \big\|_2^2  \nonumber 
\\
  \ge & 
     \frac78 \| \nabla_h \hat{e}_{\mu}^{n+1} \|_2^2   
    - 8 ( \tilde{C}_5^2 + \tilde{C}_8^2 ) ( \| e_\phi^* \|_2^2 + \| \nabla_h e_\phi^* \|_2^2 ) \quad 
  \mbox{(by~\eqref{nonlinear est-0-2}, \eqref{nonlinear est-0-4})};
\nonumber
\\
  & 
  \langle  \hat{e}_{\mu}^{n+1} , \mathbf{G}_\phi^{n+1}  \rangle_c  
  \le \| \nabla_h \hat{e}_{\mu}^{n+1} \|_2 \cdot \| \mathbf{G}_\phi^{n+1} \|_{-1,h} 
  \le \frac{1}{16} \| \nabla_h \hat{e}_{\mu}^{n+1} \|_2^2 + 4 \| \mathbf{G}_\phi^{n+1} \|_{-1,h}^2 , 
  \label{convergence-CHNS-SAV-ZEC-2-6} 
\end{align} 
in which the Cauchy inequality, as well as the preliminary nonlinear error estimates in Proposition~\ref{prop: nonlinear est}, have been repeated applied in the derivation. Subsequently, a substitution of~\eqref{convergence-CHNS-SAV-ZEC-2-1}-\eqref{convergence-CHNS-SAV-ZEC-2-6} into \eqref{convergence-CHNS-SAV-ZEC-1} leads to 
\begin{align} 
 & 
 \frac{\varepsilon^2}{4 \dt} \big( \| \nabla_h e_\phi^{n+1} \|_2^2 - \| \nabla_h e_\phi^n \|_2^2 
  +\big \| 2 \nabla_h ( e_\phi^{n+1} - e_\phi^n )\big\|_2^2 
   - \big\| 2 \nabla_h ( e_\phi^n - e_\phi^{n-1} ) \big\|_2^2  \big)   \label{convergence-CHNS-SAV-ZEC-3} 
   \\
  +& \frac34 \| \nabla_h \hat{e}_{\mu}^{n+1} \|_2^2   \nonumber
\\
  \le & 
   - \frac{1}{\dt} \langle \frac32 e_\phi^{n+1} - 2 e_\phi^n + \frac12 e_\phi^{n-1} ,  e_\mu^{n+1, (1)} \rangle_c 
   + e_q^{n+1} \langle  \phi^* \bu^* ,  \nabla_h e_{\mu}^{n+1} \rangle_1  \nonumber
\\
     +& 8 \tilde{C}_{11}^2 ( \| e_\phi^* \|_2^2 +  \| \eu^* \|_2^2 )
    + \tilde{C}_{12} ( \tilde{C}_4 + \tilde{C}_8 ) ( \frac12 | e_q^{n+1} |^2 
  + \| e_\phi^* \|_2 ^2 + \| \nabla_h e_\phi^* \|_2^2 )  \nonumber
\\
  +& 8 ( \tilde{C}_5^2 + \tilde{C}_8^2 ) ( \| e_\phi^* \|_2^2 + \| \nabla_h e_\phi^* \|_2^2 ) 
  + 4 \| \mathbf{G}_\phi^{n+1} \|_{-1,h}^2  . \nonumber
\end{align} 

In terms of the error evolutionary equation~\eqref{error-CHNS-SAV-ZEC-3}, its product with $2 e_r^{n+1}$ gives 
\begin{align}  
  & 
  \frac{1}{\dt} ( 3 e_r^{n+1} - 4 e_r^n + e_r^{n-1} ) e_r^{n+1}  
  = \frac{e_r^{n+1}}{\sqrt{E_{1,h} (\phi^*) }}  
  \Big\langle  ( \Phi_N^* )^3 - \Phi_N^* , 
   \frac{\frac32 e_\phi^{n+1} - 2 e_\phi^n + \frac12 e_\phi^{n-1}}{\dt} \Big\rangle_c    \label{convergence-CHNS-SAV-ZEC-4} 
\\
  +& \frac{e_r^{n+1}}{\sqrt{E_{1,h} (\phi^*) }}  
  \Big\langle ( {\cal NLP}^* -1 ) e_\phi^* , 
   \frac{\frac32 \Phi_N^{n+1} - 2 \Phi_N^n + \frac12 \Phi_N^{n-1}}{\dt} \Big\rangle_c 
   + 2 e_r^{n+1} \cdot \mathbf{G}_r^{n+1} \nonumber 
\\
  +&  \frac{e_r^{n+1} ( \sqrt{E_{1, h} (\phi^*)} - \sqrt{E_{1,h} (\Phi_N^*)} )}{\sqrt{E_{1, h} (\phi^*) } 
  \cdot \sqrt{E_{1, h} (\Phi_N^*) }} 
  \Big\langle ( \Phi_N^* )^3 - \Phi_N^* , 
   \frac{\frac32 \Phi_N^{n+1} - 2 \Phi_N^n + \frac12 \Phi_N^{n-1}}{\dt} \Big\rangle_c .  
\nonumber
\end{align}  
Similarly, the following estimates become available: 
\begin{align} 
  & 
  ( 3 e_r^{n+1} - 4 e_r^n + e_r^{n-1} ) e_r^{n+1}    \label{convergence-CHNS-SAV-ZEC-5-1}  
\\
  \ge & \frac12 ( | e_r^{n+1} |^2 - | e_r^n |^2 
  + | 2  e_r^{n+1} - e_r^n |^2 - | 2 e_r^n - e_r^{n-1} |^2  ) ;  \nonumber    
\\
  & 
  \frac{e_r^{n+1}}{\sqrt{E_{1,h} (\phi^*) }}  
  \Big\langle  ( \Phi_N^* )^3 - \Phi_N^* , 
   \frac{\frac32 e_\phi^{n+1} - 2 e_\phi^n + \frac12 e_\phi^{n-1}}{\dt} \Big\rangle_c \label{convergence-CHNS-SAV-ZEC-5-2}\\
   =& \frac{1}{\dt} \left\langle \frac32 e_\phi^{n+1} - 2 e_\phi^n + \frac12 e_\phi^{n-1} ,  
   e_\mu^{n+1, (1)} \right\rangle_c;
     \nonumber
\\
  & 
  \frac{e_r^{n+1}}{\sqrt{E_{1,h} (\phi^*) }}  
  \Big\langle ( {\cal NLP}^* -1 ) e_\phi^* , 
   \frac{\frac32 \Phi_N^{n+1} - 2 \Phi_N^n + \frac12 \Phi_N^{n-1}}{\dt} \Big\rangle_c \label{convergence-CHNS-SAV-ZEC-5-3} 
\\
  \le & 
    \tilde{C}_9 | e_r^{n+1} | \cdot \| e_\phi^* \|_2  
  \le \frac{\tilde{C}_9}{2} ( | e_r^{n+1} |^2 + \| e_\phi^* \|_2^2 ) \quad 
  \mbox{(by~\eqref{nonlinear est-0-5})};
  \nonumber
\\
  &
  \frac{e_r^{n+1} \bigl( \sqrt{E_{1, h} (\phi^*)} - \sqrt{E_{1,h} (\Phi_N^*)} \bigr)}{\sqrt{E_{1, h} (\phi^*) } 
  \cdot \sqrt{E_{1, h} (\Phi_N^*) }} 
  \Big\langle ( \Phi_N^* )^3 - \Phi_N^* , 
   \frac{\frac32 \Phi_N^{n+1} - 2 \Phi_N^n + \frac12 \Phi_N^{n-1}}{\dt} \Big\rangle_c \label{convergence-CHNS-SAV-ZEC-5-4}
\\
  \le &  
   \tilde{C}_{10} | e_r^{n+1} | \cdot \| e_\phi^* \|_2 
  \le \frac{\tilde{C}_{10}}{2} ( | e_r^{n+1} |^2 + \| e_\phi^* \|_2^2 ) \quad 
  \mbox{(by~\eqref{nonlinear est-0-6})};     \nonumber 
\\
  & 
  2 e_r^{n+1} \cdot \mathbf{G}_r^{n+1} 
  \le | e_r^{n+1} |^2 + | \mathbf{G}_r^{n+1} |^2 . \label{convergence-CHNS-SAV-ZEC-5-5} 
\end{align} 
In turn, a combination of \eqref{convergence-CHNS-SAV-ZEC-3}-\eqref{convergence-CHNS-SAV-ZEC-5-5} results in 
\begin{equation} 
\begin{aligned} 
 & 
 \frac{\varepsilon^2}{4 \dt} ( \| \nabla_h e_\phi^{n+1} \|_2^2 - \| \nabla_h e_\phi^n \|_2^2 
  + \| 2 \nabla_h ( e_\phi^{n+1} - e_\phi^n )\|_2^2 
   - \| 2 \nabla_h ( e_\phi^n - e_\phi^{n-1} ) \|_2^2  )  
\\
  +& \frac34 \| \nabla_h \hat{e}_{\mu}^{n+1} \|_2^2  + \frac{1}{2 \dt} ( | e_r^{n+1} |^2 - | e_r^n |^2 + | 2  e_r^{n+1} - e_r^n |^2 - | 2 e_r^n - e_r^{n-1} |^2  )
\\
  \le & 
    e_q^{n+1} \langle  \phi^* \bu^* ,  \nabla_h e_{\mu}^{n+1} \rangle_1   
   + \tilde{C}_{15} ( \| e_\phi^* \|_2 ^2 + \| \nabla_h e_\phi^* \|_2^2 ) 
   + 8 \tilde{C}_{11}^2   \| \eu^* \|_2^2  
    \\
   +& \frac12 \tilde{C}_{12} ( \tilde{C}_4 + \tilde{C}_8 ) | e_q^{n+1} |^2 
 + \frac12  (\tilde{C}_9 + \tilde{C}_{10} +2 ) | e_r^{n+1} |^2 
  + 4 \| \mathbf{G}_\phi^{n+1} \|_{-1,h}^2  + | \mathbf{G}_r^{n+1} |^2 , 
\end{aligned} 
  \label{convergence-CHNS-SAV-ZEC-6} 
\end{equation} 
with $\tilde{C}_{15} = 8 \tilde{C}_{11}^2 + \tilde{C}_{12} ( \tilde{C}_5 + \tilde{C}_8 ) + 8 ( \tilde{C}_5^2 + \tilde{C}_8^2 ) + \frac12 ( \tilde{C}_9 + \tilde{C}_{10})$.

Taking a discrete inner product with~\eqref{error-CHNS-SAV-ZEC-4} by $\ehu^{n+1}$ gives 
\begin{equation} 
\begin{aligned} 
  & 
  \frac{1}{\dt}\left \langle \frac32 \ehu^{n+1} - 2 \eu^n + \frac12 \eu^{n-1} , \ehu^{n+1}\right \rangle_1  
  + \left\langle \bu^* \cdot \nabla_h \eu^* + \eu^* \cdot \nabla_h \bU^* , 
   \ehu^{n+1} \right\rangle_1 
\\
 +& \nu \| \nabla_h \ehu^{n+1} \|_2^2 = 
    - \langle \nabla_h e_p^n , \ehu^{n+1} \rangle_1 
    - e_q^{n+1} \langle \bu^* \cdot \nabla_h \bu^* , \ehu^{n+1} \rangle_1 
\\
  +& \lambda \langle M^* \nabla_h e_\phi^* + e_{\mu}^* \nabla_h \phi^* , \ehu^{n+1} \rangle_1  + e_q^{n+1} \langle \tilde{\mu}^* \nabla_h \phi^* , \ehu^{n+1} \rangle_1 
  + \langle \mathbf{G}_{\bu}^{n+1}  ,  \ehu^{n+1} \rangle_1 . 
\end{aligned} 
  \label{convergence-CHNS-SAV-ZEC-7}  
\end{equation} 
Again, the following nonlinear error estimates could be similarly derived: 
\begin{align} 
  & 
  - \langle \eu^* \cdot \nabla_h \bU^* , 
   \ehu^{n+1} \rangle_1  \le \| \eu^* \|_2 \cdot \| \nabla_h \bU^* \|_2 
   \cdot \| \ehu^{n+1} \|_2 \le C^* \| \eu^* \|_2 \cdot \| \ehu^{n+1} \|_2 \label{convergence-CHNS-SAV-ZEC-9-1}  
\\
  \le & 
   \frac{C^*}{2} ( \| \eu^* \|_2^2 + \| \ehu^{n+1} \|_2^2 ) ;
     \nonumber 
\\
  & 
  - \langle \bu^* \cdot \nabla_h \eu^*  ,  \ehu^{n+1} \rangle_1   
  \le \| \bu^* \|_4 \cdot \| \nabla_h \eu^* \|_2  \cdot \| \ehu^{n+1} \|_4  
  \label{convergence-CHNS-SAV-ZEC-9-2}  
\\
  \le & \breve{C}_1 \tilde{C}_2 \| \nabla_h \eu^* \|_2  \cdot \| \ehu^{n+1} \|_4   
  \le  8 \breve{C}_1^2 \tilde{C}_2^2 \nu^{-1}  \| \ehu^{n+1} \|_4^2 
   + \frac{\nu}{32} \| \nabla_h \eu^* \|_2^2 ;   \nonumber
\\
  & 
  \lambda \langle M^* \nabla_h e_\phi^* + e_{\mu}^* \nabla_h \phi^* , \ehu^{n+1} \rangle_1  
  \le \lambda \| M^* \nabla_h e_\phi^* + e_{\mu}^* \nabla_h \phi^*\|_2 \cdot \| \ehu^{n+1} \|_2  \label{convergence-CHNS-SAV-ZEC-9-3}  
\\
  \le& 
     \tilde{C}_{13} \lambda ( \| \nabla_h e_\phi^* \|_2 + \| e_{\mu}^* \|_2 ) \cdot \| \ehu^{n+1} \|_2 \nonumber\\
 \le& \frac12 \tilde{C}_{13} \lambda ( \| \nabla_h e_\phi^* \|_2 ^2 
  + ( 1+  \tilde{C}_{13} \lambda) \| \ehu^{n+1} \|_2^2 ) + \frac12 \| e_{\mu}^* \|_2^2 ; 
  \nonumber 
\\
  & 
  \langle \mathbf{G}_{\bu}^{n+1}  ,  \ehu^{n+1} \rangle_1 
  \le \frac12 ( \| \mathbf{G}_{\bu}^{n+1} \|_2^2 + \|  \ehu^{n+1} \|_2^2 ) , 
  \label{convergence-CHNS-SAV-ZEC-9-4}  
\end{align} 
with an application of preliminary estimate~\eqref{nonlinear est-0-8}, as well as the functional bounds \eqref{assumption:W1-infty-1}-\eqref{assumption:W1-infty-3} and the \emph{a-priori} assumption~\eqref{a priori-1}. The temporal differentiation term could be analyzed as follow: 
\begin{equation} 
  \big\langle \frac32 \ehu^{n+1} - 2 \eu^n + \frac12 \eu^{n-1} ,  \ehu^{n+1} \big\rangle_c  
  \ge  
  \frac14 ( \| \ehu^{n+1} \|_2^2 - \| \eu^n \|_2^2 
  + \| 2 \ehu^{n+1} - \eu^n \|_2^2 
   - \| 2 \eu^n - \eu^{n-1}  \|_2^2  ) .   
    \label{convergence-CHNS-SAV-ZEC-8} 
\end{equation} 
Meanwhile, taking a discrete inner product with~\eqref{error-CHNS-SAV-ZEC-6} by $2 \eu^{n+1}$ gives 
\begin{align} 
  & 
  \| \eu^{n+1} \|_2^2 - \| \ehu^{n+1} \|_2^2 + \| \eu^{n+1} - \ehu^{n+1} \|_2^2 
  =  \| \eu^{n+1} \|_2^2 - \| \ehu^{n+1} \|_2^2   \label{convergence-CHNS-SAV-ZEC-10-1} 
\\
  +& \frac{4 \dt^2}{9} \| \nabla_h ( e_p^{n+1} - e_p^n ) \|_2^2=0,\nonumber \\
  & \mbox{so that} \quad 
   \| \ehu^{n+1} \|_2^2  =  \| \eu^{n+1} \|_2^2 
  + \frac49 \dt^2 \| \nabla_h ( e_p^{n+1} - e_p^n ) \|_2^2 , \nonumber
\end{align}  
in which the divergence-free condition for $\eu^{n+1}$ has been used. Similarly, motivated by the fact that 
\begin{equation*} 
  2 \ehu^{n+1} - \eu^n = 2 \eu^{n+1} - \eu^n + \frac43 \dt \nabla_h (p^{n+1} - p^n) , \quad 
  \nabla_h \cdot ( 2 \eu^{n+1} - \eu^n ) = 0 , 
\end{equation*} 
we are able to conclude that 
\begin{equation} 
  \| 2 \ehu^{n+1} - \eu^n \|_2^2 = \| 2 \eu^{n+1} - \eu^n \|_2^2 
  + \frac{16}{9}  \dt^2 \| \nabla_h (p^{n+1} - p^n) \|_2^2 . \label{convergence-CHNS-SAV-ZEC-10-2} 
\end{equation} 
On the other hand, in terms of the pressure gradient error term, we see that 
\begin{align} 
  \langle \nabh e_p^n , \ehu^{n+1} \rangle_1  
  = & \langle \nabla_h e_p^n , \eu^{n+1} \rangle_1 + \frac23 \dt \langle \nabla_h e_p^n , 
  \nabla_h ( e_p^{n+1} - e_p^n ) \rangle_1      \label{convergence-CHNS-SAV-ZEC-10-3} 
\\
  =  & 
  \frac23 \dt \langle \nabh e_p^n , 
  \nabla_h ( e_p^{n+1} - e_p^n ) \rangle_1    \nonumber
\\ 
  = & 
   \frac13 \dt ( \| \nabh e_p^{n+1} \|_2^2 - \| \nabla_h e_p^n \|_2^2 
  - \| \nabla_h ( e_p^{n+1} - e_p^n ) \|_2^2 )  , \nonumber
\end{align} 
in which the second step comes from the fact that $\langle \nabh e_p^n , \eu^{n+1} \rangle_1 =0$, since $\eu^{n+1}$ is divergence-free at a discrete level. Subsequently, a substitution of~\eqref{convergence-CHNS-SAV-ZEC-9-1}-\eqref{convergence-CHNS-SAV-ZEC-10-3} into \eqref{convergence-CHNS-SAV-ZEC-7} yields 
\begin{equation} 
\begin{aligned} 
  & 
  \frac{1}{4 \dt} ( \| \eu^{n+1} \|_2^2 - \| \eu^n \|_2^2 
  + \| 2 \eu^{n+1} - \eu^n \|_2^2 - \| 2 \eu^n - \eu^{n-1}  \|_2^2  ) 
  + \frac13 \dt ( \| \nabh e_p^{n+1} \|_2^2 \\
  &- \| \nabla_h e_p^n \|_2^2 )  
     + \nu \| \nabla_h \ehu^{n+1} \|_2^2  + \frac29 \dt \| \nabla_h ( e_p^{n+1} - e_p^n ) \|_2^2 
     + e_q^{n+1} \langle \bu^* \cdot \nabla_h \bu^* , \ehu^{n+1} \rangle_1 
\\
   \le&
     \frac{\nu}{32} \| \nabla_h \eu^* \|_2^2 
     + \frac12 ( C^* + \tilde{C}_{13} \lambda ( \tilde{C}_{13} \lambda + 1) + 1) \| \ehu^{n+1} \|_2^2 
     + \frac{C^*}{2} \| \eu^* \|_2^2 
\\
  &+ 8 \breve{C}_1^2 \tilde{C}_2^2 \nu^{-1}  \| \ehu^{n+1} \|_4^2 
   + \frac12 \tilde{C}_{13} \lambda  \| \nabla_h e_\phi^* \|_2 ^2 + \frac12 \| e_{\mu}^* \|_2^2  
   +  \frac12 \| \mathbf{G}_{\bu}^{n+1} \|_2^2  . 
\end{aligned} 
  \label{convergence-CHNS-SAV-ZEC-11}  
\end{equation} 

Taking a discrete inner product with~\eqref{error-CHNS-SAV-ZEC-5} by $e_q^{n+1}$ indicates that  
\begin{align} 
  & 
  \frac{1}{\dt}  \bigl( \frac32 e_q^{n+1} - 2 e_q^n + \frac12 e_q^{n-1} \bigr) e_q^{n+1}  
  + e_q^{n+1}  \langle \phi^* \bu^* , \nabla_h e_\mu^{n+1} \rangle_1   
  + e_q^{n+1}  \langle \tilde{\mu}^* \nabla_h \phi^* , \ehu^{n+1} \rangle_1  \label{convergence-CHNS-SAV-ZEC-12}  
\\
  = &
  \lambda^{-1} e_q^{n+1} \langle \bu^* \cdot \nabla_h \bu^* , \ehu^{n+1} \rangle_1  
  - e_q^{n+1} \langle e_\phi^* \bU^* + \phi^* \eu^* ,  
  \nabla_h M^{n+1} \rangle_1   + e_q^{n+1} \cdot \mathbf{G}_q^{n+1}  \nonumber 
\\
   &- e_q^{n+1} \langle M^* \nabla_h e_\phi^* + e_{\mu}^* \nabla_h \phi^* , 
   \hat{\bU}^{n+1} \rangle_1  
  + \lambda^{-1} e_q^{n+1} \langle \bu^* \cdot \nabla_h \eu^* + \eu^* \cdot \nabla_h \bU^* , 
 \hat{\bU}^{n+1} \rangle_1  .   \nonumber
\end{align} 
The following estimates could be similarly derived: 
\begin{align} 
  & 
  \bigl( \frac32 e_q^{n+1} - 2 e_q^n + \frac12 e_q^{n-1} \bigr) e_q^{n+1}   
  \label{convergence-CHNS-SAV-ZEC-13-1}  
\\
  \ge & \frac14 ( | e_q^{n+1} |^2 - | e_q^n |^2 
  + | 2  e_q^{n+1} - e_q^n |^2 - | 2 e_q^n - e_q^{n-1} |^2  ) ;  \nonumber 
\\
  & 
  - e_q^{n+1} \langle e_\phi^* \bU^* + \phi^* \eu^* ,  
  \nabla_h M^{n+1} \rangle_1  
  \le | e_q^{n+1} | \cdot \| e_\phi^* \bU^* + \phi^*\eu^*  \|_2  
  \cdot \| \nabla_h M^{n+1} \|_2  \label{convergence-CHNS-SAV-ZEC-13-2} 
\\
  \le&
   | e_q^{n+1} | \cdot \tilde{C}_{11} ( \| e_\phi^* \|_2 + \| \eu^* \|_2 ) \cdot C^{**}  \quad 
  \mbox{(by~\eqref{assumption:W1-infty-mu}, \eqref{nonlinear est-0-7})}   \nonumber 
\\ 
   \le&
     \tilde{C}_{11} C^{**} ( \frac12 | e_q^{n+1} |^2  
  + \| e_\phi^* \|_2^2 + \| \eu^* \|_2^2 );   \nonumber 
\\
  & 
  e_q^{n+1} \cdot \mathbf{G}_q^{n+1} 
  \le \frac12 ( | e_q^{n+1} |^2 + | \mathbf{G}_q^{n+1} |^2 ) , 
  \label{convergence-CHNS-SAV-ZEC-13-3} 
\\
  & 
  - e_q^{n+1} \langle M^* \nabla_h e_\phi^* + e_{\mu}^* \nabla_h \phi^* , 
   \hat{\bU}^{n+1} \rangle_1   \label{convergence-CHNS-SAV-ZEC-13-4}  
\\
   \le & | e_q^{n+1} | \cdot \| M^* \nabla_h e_\phi^* + e_{\mu}^* \nabla_h \phi^* \|_2 
   \cdot \| \hat{\bU}^{n+1} \|_2  \nonumber 
\\
  \le&
    | e_q^{n+1} | \cdot \tilde{C}_{13} ( \| \nabla_h e_\phi^* \|_2 + \| e_{\mu}^* \|_2 ) 
   \cdot | \Omega |^\frac12 \tilde{C}_2  \quad 
  \mbox{(by~\eqref{assumption:W1-infty-1}, \eqref{nonlinear est-0-8})}  \nonumber\\ 
  \le&
    \frac{1}{2} | \Omega |^\frac12 \tilde{C}_2 \tilde{C}_{13} ( \| \nabla_h e_\phi^* \|_2^2 
    + ( 1 + | \Omega |^\frac12 \tilde{C}_2 \tilde{C}_{13} ) | e_q^{n+1} |^2 )  
  +  \frac{1}{2} \| e_{\mu}^* \|_2^2 ;  \nonumber 
\\
  & 
  e_q^{n+1} \langle \bu^* \cdot \nabla_h \eu^*  ,  \hat{\bU}^{n+1} \rangle_1 
 \le\lambda^{-1}  | e_q^{n+1} | \cdot \| \bu^* \|_4 \cdot \| \nabla_h \eu^* \|_2 
 \cdot \|  \hat{\bU}^{n+1} \|_4  
 \label{convergence-CHNS-SAV-ZEC-13-5}
\\
 \le&
  \breve{C}_1 \tilde{C}_2 \cdot | \Omega |^\frac14 \tilde{C}_2 
 | e_q^{n+1} |  \cdot \| \nabla_h \eu^* \|_2  
 \le 8 | \Omega |^\frac12 \breve{C}_1^2 \tilde{C}_2^4  | e_q^{n+1} |^2 
  + \frac{\nu}{32} \| \nabla_h \eu^* \|_2^2 ;  
   \nonumber 
\\
  & 
  \lambda^{-1} e_q^{n+1} \langle \eu^* \cdot \nabla_h \bU^*  ,  \hat{\bU}^{n+1} \rangle_1 
 \le \lambda^{-1} | e_q^{n+1} | \cdot \| \eu^* \|_2 \cdot \| \nabla_h \bU^* \|_\infty 
 \cdot \|  \hat{\bU}^{n+1} \|_2  
  \label{convergence-CHNS-SAV-ZEC-13-6} 
\\
 \le&
  3 \lambda^{-1} C^* \cdot | \Omega |^\frac12 \tilde{C}_2 
 | e_q^{n+1} |  \cdot \| \eu^* \|_2  
 \le \frac{3 \lambda^{-1}}{2} | \Omega |^\frac12 C^* \tilde{C}_2  ( | e_q^{n+1} |^2 + \| \eu^* \|_2^2 )  ,    
\nonumber
\\
  & 
  - e_q^{n+1}  \langle \tilde{\mu}^* \nabla_h \phi^* , \ehu^{n+1} \rangle_1 
  \le |e_q^{n+1} | \cdot \| \tilde{\mu}^* \|_2 \cdot \| \nabla_h \phi^* \|_\infty \cdot \| \ehu^{n+1} \|_2 
   \label{convergence-CHNS-SAV-ZEC-13-7}  
\\
  \le & 
   3 \tilde{C}_2 \tilde{C}_{2,2} \cdot | e_q^{n+1} |  \cdot \| \ehu^{n+1} \|_2  
   \le \frac32 \tilde{C}_2 \tilde{C}_{2,2} ( | e_q^{n+1} |^2 + \| \ehu^{n+1} \|_2^2 ) .   
   \nonumber 
\end{align} 
In turn, a substitution of \eqref{convergence-CHNS-SAV-ZEC-13-1}-\eqref{convergence-CHNS-SAV-ZEC-13-6} into \eqref{convergence-CHNS-SAV-ZEC-12} leads to 
\begin{align} 
  & 
  \frac{1}{4 \dt} ( | e_q^{n+1} |^2 - | e_q^n |^2 
  + | 2  e_q^{n+1} - e_q^n |^2 - | 2 e_q^n - e_q^{n-1} |^2  ) 
  + e_q^{n+1}  \langle \bu^* \phi^* ,  \nabla_h e_\mu^{n+1} \rangle_c    \label{convergence-CHNS-SAV-ZEC-14}  
\\
 \le&
   \lambda^{-1} e_q^{n+1} \langle \bu^* \cdot \nabla_h \bu^* , \ehu^{n+1} \rangle_1  
  + \tilde{C}_{16} ( \| e_\phi^* \|_2^2 + \| \nabla_h e_\phi^* \|_2^2  ) 
  + \tilde{C}_{17} | e_q^{n+1} |^2  \nonumber 
\\
   &
   + \tilde{C}_{18} \| \eu^* \|_2^2  
   + \frac32 \tilde{C}_2 \tilde{C}_{2,2} \| \ehu^{n+1} \|_2^2 
   + \frac{\lambda^{-1} \nu}{32} \| \nabla_h \eu^* \|_2^2  
   + \frac{1}{2} \| e_{\mu}^* \|_2^2 
   + \frac12 | \mathbf{G}_q^{n+1} |^2 ,   \nonumber
\end{align} 
with $\tilde{C}_{16} = \tilde{C}_{11} C^{**} + \frac12 | \Omega |^\frac12 \tilde{C}_2 \tilde{C}_{13}$, $\tilde{C}_{17} = \frac12 \bigl( 1 + \tilde{C}_{11} C^{**} + 3 \tilde{C}_2 \tilde{C}_{2,2} + | \Omega |^\frac12 ( \tilde{C}_2 \tilde{C}_{13}  ( 1 + | \Omega |^\frac12 \tilde{C}_2 \tilde{C}_{13} ) + 16 \lambda^{-2} \breve{C}_1^2 \tilde{C}_2^4 + 3 \lambda^{-1} C^* \tilde{C}_2 ) \bigr)$, $\tilde{C}_{18} = \tilde{C}_{11} C^{**} + \frac{3 \lambda^{-1}}{2} | \Omega |^\frac12 C^* \tilde{C}_2$. 

Therefore, a combination of~\eqref{convergence-CHNS-SAV-ZEC-6}, \eqref{convergence-CHNS-SAV-ZEC-11} and \eqref{convergence-CHNS-SAV-ZEC-14} reveals that 
\begin{align} 
 & 
 \frac{\varepsilon^2}{4 \dt} ( \| \nabla_h e_\phi^{n+1} \|_2^2 - \| \nabla_h e_\phi^n \|_2^2 
  + \| 2 \nabla_h ( e_\phi^{n+1} - e_\phi^n )\|_2^2 
   - \| 2 \nabla_h ( e_\phi^n - e_\phi^{n-1} ) \|_2^2  )    \label{convergence-CHNS-SAV-ZEC-15} 
\\
  &+ \frac34 \| \nabla_h \hat{e}_{\mu}^{n+1} \|_2^2  + \frac{1}{2 \dt} ( | e_r^{n+1} |^2 - | e_r^n |^2 + | 2  e_r^{n+1} - e_r^n |^2 
  - | 2 e_r^n - e_r^{n-1} |^2  )  \nonumber\\
  &+ \frac{2 \lambda^{-1}}{9} \dt \| \nabla_h ( e_p^{n+1} - e_p^n ) \|_2^2 
  + \frac{\lambda^{-1}}{4 \dt} ( \| \eu^{n+1} \|_2^2 - \| \eu^n \|_2^2 
  + \| 2 \eu^{n+1} - \eu^n \|_2^2 \nonumber\\
  &- \| 2 \eu^n - \eu^{n-1}  \|_2^2  )   
  + \frac{\nu}{\lambda} \| \nabla_h \ehu^{n+1} \|_2^2  
  + \frac{\lambda^{-1}}{3} \dt ( \| \nabh e_p^{n+1} \|_2^2 - \| \nabla_h e_p^n \|_2^2 )  \nonumber \\
  &+ \frac{1}{4 \dt} ( | e_q^{n+1} |^2 - | e_q^n |^2 
  + | 2  e_q^{n+1} - e_q^n |^2 - | 2 e_q^n - e_q^{n-1} |^2  ) \nonumber \\
  \le&
   \frac{\lambda^{-1} \nu}{16} \| \nabla_h \eu^* \|_2^2 
   + \tilde{C}_{19} ( \| e_\phi^* \|_2 ^2 + \| \nabla_h e_\phi^* \|_2^2 ) 
   + \tilde{C}_{20}   | e_r^{n+1} |^2 
   + \tilde{C}_{21} \| \eu^* \|_2^2 \nonumber
\\
  &+ \tilde{C}_{22} \| \ehu^{n+1} \|_2^2  + \frac{1 + \lambda^{-1}}{2} \| e_{\mu}^* \|_2^2 
   + \tilde{C}_{23} \| \ehu^{n+1} \|_4^2 
   + \tilde{C}_{24} | e_q^{n+1} |^2  \nonumber \\
   &
   + 4 \| \mathbf{G}_\phi^{n+1} \|_{-1,h}^2 
   + \frac12 ( \lambda^{-1} \| \mathbf{G}_{\bu}^{n+1} \|_2^2 
   + | \mathbf{G}_q^{n+1} |^2 ) + | \mathbf{G}_r^{n+1} |^2 ,\nonumber 
\end{align}  
with $\tilde{C}_{19} = \tilde{C}_{15} + \frac12 \tilde{C}_{13} + \tilde{C}_{16}$, $\tilde{C}_{20} = \frac12  (\tilde{C}_9 + \tilde{C}_{10} +2 )$, $\tilde{C}_{21} = 8 \tilde{C}_{11}^2 + \frac{\lambda^{-1} C^*}{2} + \tilde{C}_{18}$, $\tilde{C}_{22} = \frac{\lambda^{-1}}{2} (C^* + \tilde{C}_{13} \lambda ( \tilde{C}_{13} \lambda + 1) +1 ) + \frac32 \tilde{C}_2 \tilde{C}_{2,2}$, $\tilde{C}_{23} = 8 \breve{C}_1^2 \tilde{C}_2^2 \lambda^{-1} \nu^{-1}$, $\tilde{C}_{24} = \frac12 \tilde{C}_{12} ( \tilde{C}_4 + \tilde{C}_8 ) + \tilde{C}_{17}$. Notice that the coupled terms have cancelled with each other, and this subtle fact has played an essential role in the theoretical proof. Meanwhile, the following estimates are valid for the star profiles: 
\begin{align} 
  \| e_\phi^* \|_2 &\le C_0 \| \nabla_h e_\phi^* \|_2 \, \, \, \mbox{(since $\overline{e_\phi^*} =0$)} ,  \label{convergence-CHNS-SAV-ZEC-16-1} \\
  &\mbox{so that} \, \, \, \| e_\phi^* \|_2 ^2 + \| \nabla_h e_\phi^* \|_2^2 
  \le ( 1 + C_0^2 ) \| \nabla_h e_\phi^* \|_2^2 , 
\nonumber \\
   \| \nabla_h e_\phi^* \|_2^2 &=  \| \nabla_h ( 2 e_\phi^n - e_\phi^{n-1} ) \|_2^2 
   \le  6 \| \nabla_h e_\phi^n \|_2^2 + 3 \| \nabla_h e_\phi^{n-1} \|_2^2 , \nonumber \\
   \| \nabla_h \Delta_h e_\phi^* \|_2^2 
  & \le  6 \| \nabla_h \Delta_h e_\phi^n \|_2^2 + 3 \| \nabla_h \Delta_h e_\phi^{n-1} \|_2^2 ,  
\nonumber \\ 
  \| \eu^* \|_2^2 &\le 6 \| \eu^n \|_2^2 + 3 \| \eu^{n-1} \|_2^2 ,  \, \, 
  \| \nabla_h \eu^* \|_2^2 \le 6 \| \nabla_h \eu^n \|_2^2 + 3 \| \nabla_h \eu^{n-1} \|_2^2 ,\nonumber
\end{align} 
where the Cauthy inequality is applied in the process. Regarding the term $\| e_{\mu}^* \|_2^2$, the following inequalities are straightforward: 
\begin{equation}
\begin{aligned} 
  \| \Delta_h e_\phi^* \|_2^2 = & - \langle \nabla_h e_\phi^* , \nabla_h \Delta_h e_\phi^* \rangle_1 
  \le 16 \| \nabla_h e_\phi^* \|_2^2 + \frac{1}{64} \| \nabla_h \Delta_h e_\phi^* \|_2^2 , 
\\
  \| e_{\mu}^* \|_2^2 \le & \tilde{C}_{14} \| e_\phi^* \|_2^2 
   + 2 \varepsilon^4 \| \Delta_h e_\phi^* \|_2^2
   \le  \tilde{C}_{14} \| e_\phi^* \|_2^2 + 16 \varepsilon^4 \| \nabla_h e_\phi^* \|_2^2 
   + \frac{\varepsilon^4}{32} \| \nabla_h \Delta_h e_\phi^* \|_2^2 
\\
  \le & 
  ( \tilde{C}_{14} + 1 ) ( \| e_\phi^* \|_2^2 + \| \nabla_h e_\phi^* \|_2^2 ) 
  + \frac{\varepsilon^4}{32} ( 6 \| \nabla_h \Delta_h e_\phi^n \|_2^2 
  + 3 \| \nabla_h \Delta_h e_\phi^{n-1} \|_2^2 ) , 
\end{aligned} 
  \label{convergence-CHNS-SAV-ZEC-16-2} 
\end{equation} 
provided that $16 \varepsilon^4 \le 1$, in which the preliminary estimate~\eqref{nonlinear est-0-8} has been recalled. In terms of the phase field diffusion part, the following estimates are observed: 
\begin{equation} 
\begin{aligned}  
  \| \nabla_h \hat{e}_{\mu}^{n+1} \|_2^2  
  = & \| \nabla_h ( - \varepsilon^2 \Delta_h e_\phi^{n+1} + e_{\mu}^{n+1, (1)} ) \|_2^2  
  \ge \frac34 \| \varepsilon^2 \nabla_h  \Delta_h e_\phi^{n+1} \|_2^2 
  - \| \nabla_h e_{\mu}^{n+1, (1)} \|_2^2  
\\
  \ge & 
   \frac34 \varepsilon^4 \| \nabla_h  \Delta_h e_\phi^{n+1} \|_2^2 
  - \tilde{C}_7^2 | e_r^{n+1} |^2 , 
\end{aligned} 
  \label{convergence-CHNS-SAV-ZEC-16-3} 
\end{equation} 
in which the Cauchy inequality has been applied in the second step, and the preliminary estimate~\eqref{nonlinear est-0-8} has been recalled in the last step. Moreover, regarding the $\| \ehu^{n+1} \|_4^2$ term, involving a $\| \cdot \|_4$ norm, we make use of~\eqref{Sobolev-1} in Proposition~\ref{prop: Poincare} and see that 
\begin{equation} 
  \tilde{C}_{23} \| \ehu^{n+1} \|_4^2 
  \le \tilde{C}_{23} C_2^2 \| \ehu^{n+1} \|_2 \cdot \| \nabla_h \ehu^{n+1} \|_2  
  \le  \tilde{C}_{23}^2 C_2^4 \lambda \nu^{-1} \| \ehu^{n+1} \|_2^2 
  + \frac{\nu}{4 \lambda} \| \nabla_h \ehu^{n+1} \|_2^2 . 
  \label{convergence-CHNS-SAV-ZEC-16-4} 
\end{equation} 
As a result, a substitution of~\eqref{convergence-CHNS-SAV-ZEC-16-1}-\eqref{convergence-CHNS-SAV-ZEC-16-4} into \eqref{convergence-CHNS-SAV-ZEC-15} yields  
\begin{align} 
 & 
 \frac{\varepsilon^2}{4 \dt} ( \| \nabla_h e_\phi^{n+1} \|_2^2 - \| \nabla_h e_\phi^n \|_2^2 
  + \| 2 \nabla_h ( e_\phi^{n+1} - e_\phi^n )\|_2^2 
   - \| 2 \nabla_h ( e_\phi^n - e_\phi^{n-1} ) \|_2^2  )     \label{convergence-CHNS-SAV-ZEC-17-1} 
\\ 
  &+ \frac{1}{2 \dt} ( | e_r^{n+1} |^2 - | e_r^n |^2 + | 2  e_r^{n+1} - e_r^n |^2 
  - | 2 e_r^n - e_r^{n-1} |^2  ) 
  + \frac{\dt}{9 \lambda} \| \nabla_h ( e_p^{n+1} - e_p^n ) \|_2^2 \nonumber
\\
  &+ \frac{\lambda^{-1}}{4 \dt} ( \| \eu^{n+1} \|_2^2 - \| \eu^n \|_2^2 
  + \| 2 \eu^{n+1} - \eu^n \|_2^2 - \| 2 \eu^n - \eu^{n-1}  \|_2^2  )  \nonumber
\\ 
  &+ \frac{\lambda^{-1}}{3} \dt ( \| \nabh e_p^{n+1} \|_2^2 - \| \nabla_h e_p^n \|_2^2 )  
  + \frac{1}{4 \dt} \bigl( | e_q^{n+1} |^2 - | e_q^n |^2 
  + | 2  e_q^{n+1} - e_q^n |^2 \nonumber \\
  &- | 2 e_q^n - e_q^{n-1} |^2  \bigr)  
  +\frac{3 \nu}{4 \lambda} \| \nabla_h \eu^{n+1} \|_2^2  
  - \frac{3 \nu}{8 \lambda}  \| \nabla_h \eu^n \|_2^2 
  - \frac{\nu}{16 \lambda}  \| \nabla_h \eu^{n-1} \|_2^2   \nonumber
  \\
  & 
  + \frac{9}{16} \varepsilon^4 \| \nabla_h  \Delta_h e_\phi^{n+1} \|_2^2 
  - \frac{\varepsilon^4}{32} ( 6 \| \nabla_h \Delta_h e_\phi^n \|_2^2 
  + 3 \| \nabla_h \Delta_h e_\phi^{n-1} \|_2^2 )  \nonumber \\
  \le&
     \tilde{C}_{25} ( 6 \| \nabla_h e_\phi^n \|_2^2 + 3 \| \nabla_h e_\phi^{n-1} \|_2^2 )  + \tilde{C}_{24} | e_q^{n+1} |^2+ \tilde{C}_{21} ( 6 \| \eu^n \|_2^2 + 3 \| \eu^{n-1} \|_2^2 ) 
    \nonumber \\
   &
   + \tilde{C}_{27} \| \eu^{n+1} \|_2^2  + \tilde{C}_{26}   | e_r^{n+1} |^2 
   + 4 \| \mathbf{G}_\phi^{n+1} \|_{-1,h}^2  \nonumber 
\\
  & + \frac12 ( \lambda^{-1} \| \mathbf{G}_{\bu}^{n+1} \|_2^2 
   + | \mathbf{G}_q^{n+1} |^2 ) 
   + | \mathbf{G}_r^{n+1} |^2 , \nonumber
\end{align}  
with $\tilde{C}_{25} = (\tilde{C}_{19} + \tilde{C}_{14} + 1) ( 1 + C_0^2)$, $\tilde{C}_{26} = \tilde{C}_{20} + \frac34 \tilde{C}_7^2$, $\tilde{C}_{27} = \tilde{C}_{22} + \tilde{C}_{23}^2 C_2^4 \lambda \nu^{-1}$. In fact, the inequality $ \| \nabla_h \eu^{n+1} \|_2^2 \le \| \nabla_h \ehu^{n+1} \|_2^2$ (as indicated by~\eqref{projection est-0-1} in Proposition~\ref{prop: projection est}), the identity $\| \ehu^{n+1} \|_2^2  =  \| \eu^{n+1} \|_2^2 
  + \frac49 \dt^2 \| \nabla_h ( e_p^{n+1} - e_p^n ) \|_2^2$ by~\eqref{convergence-CHNS-SAV-ZEC-10-1}, as well as the bound that $\tilde{C}_{27} \cdot \frac49 \dt^2 \le \frac{\lambda^{-1}}{9} \dt$ (provided that $\dt$ is sufficiently small), have been applied in the derivation. 
  
The following quantity is introduced for the convenience of the convergence analysis: 
\begin{equation} 
\begin{aligned} 
 F^{n+1} &:=  \frac{\varepsilon^2}{4} ( \| \nabla_h e_\phi^{n+1} \|_2^2  
  + \| 2 \nabla_h ( e_\phi^{n+1} - e_\phi^n )\|_2^2 )    
  + \frac{\lambda^{-1}}{4} ( \| \eu^{n+1} \|_2^2 + \| 2 \eu^{n+1} - \eu^n \|_2^2 )    
\\
  &   
  + \frac{\lambda^{-1}}{3} \dt^2 \| \nabh e_p^{n+1} \|_2^2 
  + \frac12 ( | e_r^{n+1} |^2 + | 2  e_r^{n+1} - e_r^n |^2  )  
\\
  & 
  + \frac14 ( | e_q^{n+1} |^2 + | 2  e_q^{n+1} - e_q^n |^2   )  . 
\end{aligned} 
  \label{convergence-CHNS-SAV-ZEC-17-2} 
\end{equation} 
Going back \eqref{convergence-CHNS-SAV-ZEC-17-1}, we see that 
\begin{equation} 
\begin{aligned} 
 & 
 \frac{1}{\dt} ( F^{n+1} - F^n ) 
  + \frac{9}{16} \varepsilon^4 \| \nabla_h \Delta_h e_\phi^{n+1} \|_2^2   
  - \frac{\varepsilon^4}{32} ( 6 \| \nabla_h \Delta_h e_\phi^n \|_2^2 
  + 3 \| \nabla_h \Delta_h e_\phi^{n-1} \|_2^2 ) 
\\
  & + \frac{3 \nu}{4 \lambda} \| \nabla_h \eu^{n+1} \|_2^2  
  - \frac{3 \nu}{8 \lambda}  \| \nabla_h \eu^n \|_2^2 
  - \frac{\nu}{16 \lambda}  \| \nabla_h \eu^{n-1} \|_2^2  
  + \frac{\lambda^{-1}}{9} \dt \| \nabla_h ( e_p^{n+1} - e_p^n ) \|_2^2  
\\
  \le & 
    A_1 F^{n+1} + A_0 F^n + A_{-1} F^{n-1}  
   + 4 \| \mathbf{G}_\phi^{n+1} \|_{-1,h}^2  
\\
  & 
   + \frac12 ( \lambda^{-1} \| \mathbf{G}_{\bu}^{n+1} \|_2^2 
   + | \mathbf{G}_q^{n+1} |^2 ) + | \mathbf{G}_r^{n+1} |^2 , 
\end{aligned} 
  \label{convergence-CHNS-SAV-ZEC-17-3} 
\end{equation} 
with $A_1 = \max( 4 \tilde{C}_{24} , 2 \tilde{C}_{26} , 4 \tilde{C}_{27} \lambda )$, $A_0 = \max ( 24 \tilde{C}_{25} \varepsilon^{-2} , 24 \tilde{C}_{21} )$, and $A_{-1} = \\ \max ( 12  \tilde{C}_{25} \varepsilon^{-2} , 12 \tilde{C}_{21} \lambda)$. Again, all these constants only depends on the regularity of the exact solution, the domain $\Omega$ and the physical parameters. In turn, an application of discrete Gronwall inequality leads to the desired error estimate  
\begin{align} 
  & 
  F^{n+1} + \frac{\varepsilon^4}{4} \dt \sum_{k=1}^{n+1}  \| \nabla_h \Delta_h e_\phi^k \|_2^2 
  + \frac{\nu}{4 \lambda} \dt \sum_{k=1}^{n+1}  \| \nabla_h \eu^k \|_2^2  
   \le \hat{C}_1 ( \dt^2 + h^2 )^2 ,  \, \,  
  \mbox{so that}\label{convergence-CHNS-SAV-ZEC-17-4} 
\\
  & 
  \varepsilon \| \nabla_h \ephi^{n+1} \|_2 + \| \eu^{n+1} \|_2 
  + \varepsilon^2 \Big( \dt \sum_{k=1}^{n+1}  \| \nabla_h \Delta_h e_\phi^k \|_2^2 \Big)^\frac12  
  + \Big( \frac{\nu}{\lambda} \dt \sum_{k=1}^{n+1}  \| \nabla_h \eu^k \|_2^2 \Big)^\frac12  
\nonumber\\
   \le &2 \hat{C}_1^\frac12 ( \dt^2 + h^2 ) , \nonumber
\end{align}  
in which the accuracy order of the local truncation errors has been used. As a result, an optimal rate error estimate is obtained.  

\subsection{Recovery of the \emph{a-priori} assumption~\eqref{a priori-1} at the next time step} 

With the full order convergence estimate \eqref{convergence-CHNS-SAV-ZEC-17-4} in hand, the \emph{a-priori} assumption in~\eqref{a priori-1} could be appropriately recovered. The analysis is separately performed in two different cases, in terms of the scaling law between the time step and spatial mesh sizes. 

If $\dt \le h$, an application of inverse inequality implies that 
\begin{equation} 
\begin{aligned} 
  \| \eu^{n+1} \|_4& \le  \frac{C \| \eu^{n+1} \|_2}{h^\frac12} 
  \le   \frac{2 C \hat{C}_1 (\dt^2 + h^2) }{h^\frac12}  \le \dt^\frac14 + h^\frac14 , 
\\
  \| \nabla_h e_\phi^{n+1} \|_2 &\le 2 \hat{C}_1^\frac12 \varepsilon^{-1} ( \dt^2 + h^2) ,\\
  \| e_\phi^{n+1} \|_2 &\le C_0 \| \nabla_h e_\phi^{n+1} \|_2 
  \le 2 C_0 \hat{C}_1^\frac12 \varepsilon^{-1} ( \dt^2 + h^2) , 
\\ 
  \| e_\phi^{n+1} \|_\infty + \| \nabla_h e_\phi^{n+1} \|_\infty 
  &\le \frac{C ( \| e_\phi^{n+1} \|_2 + \| \nabla_h e_\phi^{n+1} \|_2 ) }{h}  
\\
  \le & \frac{2 ( 1 + C_0 ) \hat{C}_1^\frac12 \varepsilon^{-1} ( \dt^2 + h^2) }{h}  
  \le \dt^\frac14 + h^\frac14 , 
\end{aligned} 
	\label{a priori-4}  
\end{equation} 
provided that $\dt$ and $h$ are sufficiently small. This has validated the \emph{a-priori} assumption~\eqref{a priori-1} if $\dt \le h$. 

Conversely, if $\dt \ge h$, the diffusion error estimate in \eqref{convergence-CHNS-SAV-ZEC-17-4} reveals that 
\begin{equation} 
  \| \nabla_h \Delta_h e_\phi^{n+1} \|_2  + \| \nabla_h \eu^{n+1} \|_2  
   \le \frac{2 \hat{C}_1^\frac12 ( \varepsilon^{-2} + \lambda^\frac12 \nu^{-\frac12} ) ( \dt^2 + h^2 ) }{\dt^\frac12} \le \dt + h .  
	\label{a priori-5}
\end{equation} 
An application of Sobolev interpolation inequalities~\eqref{Sobolev-0} and \eqref{Sobolev-1} (in Proposition~\ref{prop: Poincare}) leads to 
\begin{align} 
  \| \eu^{n+1} \|_4 &\le C_2 \| \eu^{n+1} \|_2^\frac12 \cdot \| \nabla_h \eu^{n+1} \|_2 
  \le 2^\frac12 \hat{C}_1^\frac14 C_2 (\dt^2 + h^2)^\frac12 \cdot ( \dt + h )^\frac12   \label{a priori-6} \\
  &\le \dt^\frac14 + h^\frac14 , \nonumber
\\
  \| e_\phi^{n+1} \|_\infty &\le C_1 \| e_\phi^{n+1} \|_2^\frac23 \cdot \| \nabla_h \Delta_h e_\phi^{n+1} \|_2^\frac13\nonumber \\ 
  &\le 2^\frac23 C_0^\frac23 \hat{C}_1^\frac13 \varepsilon^{-\frac23} C_1 (\dt^2 + h^2)^\frac12 
  \cdot ( \dt + h )^\frac13 
   \le \frac12 ( \dt^\frac14 + h^\frac14 ) , \nonumber
\\ 
  \| \nabla_h e_\phi^{n+1} \|_\infty &\le C_1 \| \nabla_h e_\phi^{n+1} \|_2^\frac12 \cdot \| \nabla_h \Delta_h e_\phi^{n+1} \|_2^\frac12\nonumber \\
  &\le 2^\frac12 \hat{C}_1^\frac12 \varepsilon^{-\frac12} C_1 (\dt^2 + h^2)^\frac12 
  \cdot ( \dt + h )^\frac12 
   \le \frac12 ( \dt^\frac14 + h^\frac14 ) , \nonumber
\end{align}  
provided that $\dt$ and $h$ are sufficiently small, so that the \emph{a-priori} assumption~\eqref{a priori-1} has also been validated if $\dt \ge h$. As a result, an induction analysis could be effectively applied, and proof of Theorem~\ref{thm: convergence} is finished.

\section{Conclusion}
In this paper, we have rigorously derived error estimates for a fully discrete, unconditionally energy-stable scheme for the Cahn-Hilliard-Navier-Stokes system, a phase-field model for two-phase incompressible flow. The numerical is based on the scalar auxiliary variable reformulation, combined with the zero energy contribution approach. Because of this reformulation, all the nonlinear and coupled terms could be explicitly computed in the resulting numerical scheme, and only constant coefficient Poisson solvers are needed in the numerical implementation. We have established a second-order convergence rate for the proposed numerical scheme in both time and space, in the $\ell^{\infty}(0,T;H_h^1)\cap\ell^2(0,T;H_h^3)$ norm for the phase variable and the $\ell^{\infty}(0,T;\ell^2)\cap\ell^2(0,T;H_h^1)$ norm for the velocity variable, following the energy norms of the reformulated PDE system. This is the first work to establish an optimal convergence estimate for the Cahn-Hilliard-Navier-Stokes system using a ZEC-based fully decoupled scheme.

\appendix

\section{Proof of Proposition~\ref{prop: nonlinear est}} 
\label{appA}

Based on the representation formula for ${\cal NLP}^*$ (in~\eqref{error-CHNS-SAV-ZEC-2}), its $\| \cdot \|_\infty$ following bound is obvious: 
\begin{equation} 
  \| {\cal NLP}^* \|_\infty \le \frac32 ( \| \Phi_N^* \|_\infty^2 + \| \phi^* \|_\infty^2 )  
  \le \frac{27}{2} ( (C^*)^2 + \tilde{C}_2^2 ) ,  \label{nonlinear est-1-1}
\end{equation} 
with the help of~\eqref{assumption:W1-infty-2} and \eqref{a priori-3}. In terms of its gradient estimate, we see that the difference approximation expansion implies that 
\begin{equation} 
\begin{aligned} 
  \| \nabla_h {\cal NLP}^* \|_\infty \le & 
  2 ( \| \Phi_N^* \|_\infty \cdot   \| \nabla_h \Phi_N^* \|_\infty 
  + \| \phi^* \|_\infty \cdot   \| \nabla_h \phi^* \|_\infty  ) 
\\
  & 
  +  \| \Phi_N^* \|_\infty \cdot   \| \nabla_h \phi^* \|_\infty 
  + \| \phi^* \|_\infty \cdot   \| \nabla_h \Phi_N^* \|_\infty  \\
 \le &18 ( (C^*)^2 + \tilde{C}_2^2  + C^* \tilde{C}_2 ) . 
\end{aligned} 
  \label{nonlinear est-1-2}
\end{equation} 
Therefore, a combination of \eqref{nonlinear est-1-1} and \eqref{nonlinear est-1-2} yields the first inequality in \eqref{nonlinear est-0-1}, by taking $\tilde{C}_3 = \frac{27}{2} ( (C^*)^2 + \tilde{C}_2^2 ) + 18 ( (C^*)^2 + \tilde{C}_2^2  + C^* \tilde{C}_2 )$. 

The proof of the second inequality in \eqref{nonlinear est-0-1} is similar, based on the following bounds: 
\begin{equation} 
\begin{aligned} 
  & 
  \| ( \Phi_N^* )^3 \|_\infty \le \| \Phi_N^* \|_\infty^3 \le 27 ( C^*)^3  , \quad 
  \| \Phi_N^* \|_\infty \le 3 C^* ,   \quad\| \nabla_h \Phi_N^* \|_\infty \le  3 C^* ,
\\
  & 
  \| \nabla_h ( ( \Phi_N^* )^3 ) \|_\infty \le 3 \| \Phi_N^* \|_\infty^2 \cdot \| \nabla \Phi_n^* \|_\infty 
   \le 3 \cdot ( 3 C^* )^2 \cdot 3 C^* = 81 (C^*)^3. 
\end{aligned} 
  \label{nonlinear est-1-3}
\end{equation}    
This gives the the second inequality in \eqref{nonlinear est-0-1}, by taking $\tilde{C}_4 = 108 ( C^*)^3  + 6 C^*$. 

In terms of inequality \eqref{nonlinear est-0-2}, an application of discrete H\"older inequality indicates that 
\begin{equation} 
\begin{aligned} 
  \| e_\mu^{n+1, (2)} \|_2 \le & \frac{ | r^{n+1} |}{\sqrt{E_{1, h} (\phi^*) }} 
  \cdot ( \| {\cal NLP}^* \|_\infty +1 ) \| e_\phi^*\|_2  
  \le \sqrt{2} \tilde{C}_0^\frac12 \cdot | \Omega |^{-\frac12} ( \tilde{C}_3 +1 ) \| e_\phi^*\|_2  , 
\end{aligned} 
  \label{nonlinear est-2-1}
\end{equation}
in which the preliminary estimates~\eqref{energy stab-1}, \eqref{nonlinear est-0-1}, as well as the fact that $\sqrt{E_{1, h} (\phi^*)} \ge | \Omega |$, have been applied. Moreover, because of the fact that both $r^{n+1}$ and $\sqrt{E_{1, h} (\phi^*) }$ are scalar constants, the gradient estimate could be derived as follows 
\begin{equation} 
\begin{aligned} 
  \| \nabla_h e_\mu^{n+1, (2)} \|_2 \le & \frac{ | r^{n+1} |}{\sqrt{E_{1, h} (\phi^*) }} 
  \cdot \Big( ( \| {\cal NLP}^* \|_\infty +1 ) \| \nabla_h e_\phi^*\|_2  
  + \| \nabla_h {\cal NLP}^* \|_\infty  \cdot \| e_\phi^*\|_2 \Big) 
\\
  \le & 
  \sqrt{2} \tilde{C}_0^\frac12 \cdot | \Omega |^{-\frac12} \Big( ( \tilde{C}_3 +1 ) \| \nabla_h e_\phi^*\|_2  
  + \tilde{C}_3 \| e_\phi^* \|_2 ) \Big) . 
\end{aligned} 
  \label{nonlinear est-2-2}
\end{equation}
Subsequently, inequality \eqref{nonlinear est-0-2} has been proved, by taking $\tilde{C}_5 = \sqrt{2} \tilde{C}_0^\frac12 \cdot | \Omega |^{-\frac12}  ( \tilde{C}_3 +1 )$. 

To establish the nonlinear energy error estimate~\eqref{nonlinear est-0-3}, we begin with the following expansion: 
\begin{equation} 
  E_{1, h} (\Phi_N^*) - E_{1,h} (\phi^*) = \frac14 \langle ( (\Phi_N^*)^3 + (\Phi_N^*)^2 \phi^* 
  + \Phi_N^* ( \phi^* )^2 + ( \phi^* )^3 - 2 ( \Phi_N^* + \phi^* ) , e_\phi^* \rangle . 
  \label{nonlinear est-3-1}
\end{equation}
In turn, an application of discrete H\"older inequality reveals that 
\begin{equation} 
\begin{aligned} 
  | E_{1, h} (\Phi_N^*) - E_{1,h} (\phi^*) | \le & 
   \frac14 ( \| \Phi_N^* \|_\infty^3 + \| \Phi_N^* \|_\infty^2 \cdot \| \phi^* \|_\infty 
  + \| \Phi_N^* \|_\infty \cdot \| \phi^* \|_\infty^2 + \| \phi^* \|_\infty^3 
\\
  & 
  + 2 ( \| \Phi_N^* \|_\infty + \| \phi^* \|_\infty ) ) \cdot \| e_\phi^* \|_1 
\\
  \le & 
  \frac14 (  (C^*)^3 + \tilde{C}_2^3 + ( C^* \tilde{C}_2 + 2 ) ( C^* + \tilde{C}_2 ) ) 
  | \Omega |^\frac12 \cdot \| e_\phi^* \|_2 , 
\end{aligned} 
  \label{nonlinear est-3-2}
\end{equation}
in which the preliminary estimates \eqref{assumption:W1-infty-2} and \eqref{a priori-3}, as well as the fact that $\| f \|_1 \le | \Omega |^\frac12 \| f \|_2$, have been applied in the derivation. Moreover, this estimate could be used to derive the second inequality in \eqref{nonlinear est-0-3}: 
\begin{equation} 
\begin{aligned} 
  & 
  | \sqrt{E_{1, h} (\Phi_N^*)} - \sqrt{E_{1,h} (\phi^*)} |  
  = \frac{ | E_{1, h} (\Phi_N^*) - E_{1,h} (\phi^*) | }{ \sqrt{E_{1, h} (\Phi_N^*)} 
   + \sqrt{E_{1,h} (\phi^*)} }  
\\ 
  \le & 
  \frac14 (  (C^*)^3 + \tilde{C}_2^3 + ( C^* \tilde{C}_2 + 2 ) ( C^* + \tilde{C}_2 ) ) 
  | \Omega |^\frac12 \| e_\phi^* \|_2 \cdot ( 2 | \Omega |^\frac12 )^{-1} ,  
\end{aligned} 
  \label{nonlinear est-3-3}
\end{equation}
in which the fact that $E_{1, h} (f) \ge | \Omega |$ (for any $f$) has been used again. As a result, a combination of \eqref{nonlinear est-3-2} and \eqref{nonlinear est-3-3} leads to the desired inequality \eqref{nonlinear est-0-3}, by taking $\tilde{C}_6 = \frac14 (  (C^*)^3 + \tilde{C}_2^3 + ( C^* \tilde{C}_2 + 2 ) ( C^* + \tilde{C}_2 ) )  \max ( | \Omega|^\frac12 , 2^{-\frac12} )$. 

  The inequalities in \eqref{nonlinear est-0-4} could be similarly proved, and the technical details are left to interested readers. 

  The proof of \eqref{nonlinear est-0-5} and \eqref{nonlinear est-0-6} is based an application of discrete H\"older inequality: 
\begin{equation} 
\begin{aligned} 
    & 
  \Big| \frac{1}{\sqrt{E_{1,h} (\phi^*) }}  
  \Big\langle ( {\cal NLP}^* -1 ) e_\phi^* , 
   \frac{\frac32 \Phi_N^{n+1} - 2 \Phi_N^n + \frac12 \Phi_N^{n-1}}{\dt} \Big\rangle_c \Big| 
\\
   \le & 
   | \Omega |^{-\frac12} (  \| {\cal NLP}^* \|_\infty + 1 ) \| e_\phi^* \|_1  
   \cdot \Big\| \frac{\frac32 \Phi_N^{n+1} - 2 \Phi_N^n + \frac12 \Phi_N^{n-1}}{\dt} \Big\|_\infty 
\\
   \le & 
   | \Omega |^{-\frac12} (  \tilde{C}_3 + 1 ) \cdot | \Omega |^\frac12 \| e_\phi^* \|_2  \cdot C^* , 
\end{aligned} 
  \label{nonlinear est-5-1} 
\end{equation} 
\begin{equation}  
\begin{aligned} 
  & 
  \Big|  \frac{\sqrt{E_{1, h} (\phi^*)} - \sqrt{E_{1,h} (\Phi_N^*)}}{\sqrt{E_{1, h} (\phi^*) } 
  \cdot \sqrt{E_{1, h} (\Phi_N^*) }} 
  \Big\langle ( \Phi_N^* )^3 - \Phi_N^* , 
   \frac{\frac32 \Phi_N^{n+1} - 2 \Phi_N^n + \frac12 \Phi_N^{n-1}}{\dt} \Big\rangle_c  \Big| 
\\
  \le & 
  | \Omega |^{-1} \cdot  | \sqrt{E_{1, h} (\phi^*)} - \sqrt{E_{1,h} (\Phi_N^*)} | 
  \cdot \| ( \Phi_N^* )^3 - \Phi_N^* \|_\infty 
   \cdot \Big\| \frac{\frac32 \Phi_N^{n+1} - 2 \Phi_N^n + \frac12 \Phi_N^{n-1}}{\dt} \Big\|_\infty 
\\
   \le & 
   | \Omega |^{-1} \cdot \tilde{C}_6 \| e_\phi^* \|_2 \cdot \tilde{C}_4  \cdot C^* , 
\end{aligned} 
  \label{nonlinear est-6-1} 
\end{equation}  
in which the preliminary estimates \eqref{assumption:W1-infty-3}, \eqref{nonlinear est-0-1}, \eqref{nonlinear est-0-3} have been repeatedly applied. In turn, inequalities  \eqref{nonlinear est-0-5} and \eqref{nonlinear est-0-6} become valid, by taking $\tilde{C}_9 = C^* (\tilde{C}_3 +1)$, $\tilde{C}_{10} = C^* \tilde{C}_4 \tilde{C}_6  | \Omega |^{-1}$. 

The derivation of the two inequalities in \eqref{nonlinear est-0-7} is more straightforward: 
\begin{equation} 
\begin{aligned} 
  & 
  \| e_\phi^* \bU^* \|_2 \le \| e_\phi^* \|_2 \cdot \| \bU^* \|_\infty  \le 3 C^* \| e_\phi^* \|_2  , \quad  
  \| \phi^* \eu^* \|_2  \le  \| \phi^* \|_\infty \cdot \| \eu^* \|_2 \le 3 \tilde{C}_2 \| \eu^* \|_2 , 
\\
  &  
  \| \bu^* \phi^* \|_2 \le \| \bu^* \|_2 \cdot \| \phi^* \|_\infty  
  \le  | \Omega |^\frac14 \| \bu^* \|_4 \cdot \| \phi^* \|_\infty 
  \le 3 \tilde{C}_2 \cdot  \breve{C}_1 \tilde{C}_2 \cdot | \Omega |^\frac14 ,  
\end{aligned} 
  \label{nonlinear est-7-1} 
\end{equation}  
with the preliminary assumption~\eqref{assumption:W1-infty-2} and the a-priori estimate~\eqref{a priori-3} repeatedly used. Consequently, the two inequalities in \eqref{nonlinear est-0-7} are proved, by taking $\tilde{C}_{11} = 3 \max ( C^*, \tilde{C}_2 )$, $\tilde{C}_{12} = 3  \breve{C}_1 \tilde{C}_2^2 | \Omega |^\frac14$. 

The first inequality in~\eqref{nonlinear est-0-8} could be proved in a similar fashion: 
\begin{equation} 
\begin{aligned} 
  & 
  \| M^* \nabla_h e_\phi^* \|_2 \le \| M^* \|_\infty \cdot \| \nabla_h e_\phi^* \|_2 
  \le C^{**} \| \nabla_h e_\phi^* \|_2  , \quad (\mbox{by \eqref{assumption:W1-infty-3}}) 
\\
  & 
    \| e_{\mu}^* \nabla_h \phi^* \|_2  \le     \| e_{\mu}^* \|_2 \cdot \| \nabla_h \phi^* \|_\infty 
    \le  3 \tilde{C}_2 \| e_{\mu}^* \|_2  , \quad (\mbox{by \eqref{a priori-3}}) , 
 \end{aligned} 
  \label{nonlinear est-8-1}  
\end{equation} 
by taking $\tilde{C}_{13} = \max ( C^{**} , 3 \tilde{C}_2 )$. In terms of the second inequality in~\eqref{nonlinear est-0-8}, we begin with the following expansion: 
\begin{equation} 
  e_{\mu}^* = ( {\cal NLP}^* -1 ) e_\phi^* - \varepsilon^2 \Delta_h e_\phi^* . 
  \label{nonlinear est-8-2}  
\end{equation} 
In turn, a careful application of Cauchy inequality and discrete H\"older inequality gives 
\begin{equation} 
\begin{aligned} 
  &
  \| ( {\cal NLP}^* -1 ) e_\phi^* \|_2 \le ( \| {\cal NLP}^* \|_\infty + 1 ) \| e_\phi^* \|_2 
  \le ( \tilde{C}_3 +1 ) \| e_\phi^* \|_2 , 
\\
  & 
   \| e_{\mu}^* \|_2^2 \le 
   2 ( \| ( {\cal NLP}^* -1 ) e_\phi^* \|_2^2 
   + \varepsilon^4 \| \Delta_h e_\phi^* \|_2^2 ) 
   \le  2 ( \tilde{C}_3 + 1)^2 \| e_\phi^* \|_2^2 
   + 2 \varepsilon^4 \| \Delta_h e_\phi^* \|_2^2 ,  
\end{aligned} 
  \label{nonlinear est-8-3}  
\end{equation} 
which is exactly the second inequality in~\eqref{nonlinear est-0-8}, by taking $\tilde{C}_{14} = 2 ( \tilde{C}_3 + 1)^2$. The proof of Proposition~\ref{prop: nonlinear est} has been completed.

\section{Proof of Proposition~\ref{prop: projection est}}  
\label{appB}
For $\ehu^{n+1}$ and $\eu^{n+1}$ satisfying \eqref{error-CHNS-SAV-ZEC-6}-\eqref{error-CHNS-SAV-ZEC-7}, we have 
\begin{equation}  
  \| \ehu^{n+1} \|_2^2 = \| \eu^{n+1} \|_2^2 + \| \ehu^{n+1} - \eu^{n+1} \|_2^2 , ~
  \| \nabla_h \ehu^{n+1} \|_2^2 = \| \nabla_h \eu^{n+1} \|_2^2 + \| \nabla_h ( \ehu^{n+1} - \eu^{n+1} ) \|_2^2 . 
\end{equation} 

Taking a discrete inner product with \eqref{error-CHNS-SAV-ZEC-6} by $2 \eu^{n+1}$ gives 
\begin{equation} 
  \langle \eu^{n+1} - \ehu^{n+1} , 2 \eu^{n+1} \rangle_1 
  + \frac43 \dt \langle \nabla_h (e_p^{n+1} - e_p^n ) ,  \eu^{n+1} \rangle_1 = 0 . 
  \label{projection est-1-1} 
\end{equation} 
The first term on the left hand side could be expanded in a standard way: 
\begin{equation} 
  \langle \eu^{n+1} - \ehu^{n+1} , 2 \eu^{n+1} \rangle_1  
  = \| \eu^{n+1} \|_2^2 - \| \ehu^{n+1} \|_2^2 + \| \ehu^{n+1} - \eu^{n+1} \|_2^2 . 
  \label{projection est-1-2} 
\end{equation} 
The second term on the left hand side disappears, due to the discrete divergence-free identity \eqref{error-CHNS-SAV-ZEC-7} for $\eu^{n+1}$, combined with the no-penetration boundary condition, $( \eu^{n+1} \cdot \n ) \mid_{\partial \Omega} = 0$; 
\begin{equation} 
   \langle \nabla_h (e_p^{n+1} - e_p^n ) ,  \eu^{n+1} \rangle_1 
   = - \langle e_p^{n+1} - e_p^n ,  \nabla_h \cdot \eu^{n+1} \rangle_1  = 0 . 
   \label{projection est-1-3} 
\end{equation} 
In turn, a combination of \eqref{projection est-1-2} and \eqref{projection est-1-3} leads to the first equality in \eqref{projection est-0-1}. 

In addition, taking a discrete inner product with \eqref{error-CHNS-SAV-ZEC-6} by $-2 \Delta_h \eu^{n+1}$ yields 
\begin{equation} 
  - 2 \langle \eu^{n+1} - \ehu^{n+1} , \Delta_h \eu^{n+1} \rangle_1 
  - \frac43 \dt \langle \nabla_h (e_p^{n+1} - e_p^n ) ,  \Delta_h \eu^{n+1} \rangle_1 = 0 . 
  \label{projection est-2-1} 
\end{equation} 
The first term could be analyzed with the help of summation-by-parts formula: 
\begin{equation} 
\begin{aligned} 
  & 
  - 2 \langle \eu^{n+1} - \ehu^{n+1} , \Delta_h \eu^{n+1} \rangle_1  
  = 2 \langle \nabla_h ( \eu^{n+1} - \ehu^{n+1} ) , \Delta_h \eu^{n+1} \rangle_1 
\\
  = &  
   \| \nabla_h \eu^{n+1} \|_2^2 - \| \nabla_h \ehu^{n+1} \|_2^2 + \| \nabla_h ( \ehu^{n+1} - \eu^{n+1} ) \|_2^2 , 
\end{aligned} 
  \label{projection est-2-2} 
\end{equation} 
in which the no penetration, free slip boundary condition, $( \eu^{n+1} \cdot \n ) \mid_{\partial \Omega} = 0$, $\partial_{\n} ( \eu^{n+1} \cdot \boldsymbol{\tau}) \mid_{\partial \Omega} = 0$, has played an important role in the derivation. 
Meanwhile, we see that the second term on the left hand side disappears, due to the fact that $\nabla_h \cdot (\Delta_h \eu^{n+1} ) = \Delta_h ( \nabla_h \cdot \eu^{n+1} ) =0$, combined with the no penetration, free slip boundary condition for $\eu^{n+1}$: 
\begin{equation} 
\begin{aligned} 
   \langle \nabla_h (e_p^{n+1} - e_p^n ) ,  \Delta_h \eu^{n+1} \rangle_1 
   = & 
   - \langle e_p^{n+1} - e_p^n ,  \nabla_h \cdot ( \Delta_h \eu^{n+1} ) \rangle_1  = 0 
 \\
   = & 
   - \langle e_p^{n+1} - e_p^n ,  \Delta_h ( \nabla_h \cdot \eu^{n+1} ) \rangle_1  = 0 . 
 \end{aligned} 
   \label{projection est-2-3} 
\end{equation} 
In fact, the no penetration, free slip boundary condition for $\eu^{Pn+1}$ ensures that the normal component of $\Delta_h \eu^{n+1}$ vanishes on the boundary, namely, $( \Delta_h \eu^{n+1} \cdot \n ) \mid_{\partial \Omega} = 0$. This subtle fact has played an essential in the derivation of \eqref{projection est-2-3}. Therefore, a combination of \eqref{projection est-2-2} and \eqref{projection est-2-3} yields the second equality in \eqref{projection est-0-1}. The proof of Proposition~\ref{prop: projection est} is finished. 

\end{document}